\documentclass[12pt,a4paper]{article}
\usepackage{amsmath}
\usepackage{amsfonts}
\usepackage{amssymb}
\usepackage{graphicx}%
\providecommand{\U}[1]{\protect\rule{.1in}{.1in}}
\newtheorem{theorem}{Theorem}[section]

\newtheorem{algorithm}[theorem]{Algorithm}

\newtheorem{corollary}[theorem]{Corollary}

\newtheorem{definition}[theorem]{Definition}

\newtheorem{lemma}[theorem]{Lemma}

\newtheorem{remark}[theorem]{Remark}

\newenvironment{proof}[1][Proof]{\textbf{#1.} }{\ \rule{0.5em}{0.5em}}

\makeatletter\@addtoreset{equation}{section}\makeatother
\evensidemargin=\oddsidemargin
\newdimen\dummy
\dummy=\oddsidemargin
\begin{document}

\title{Variable Order, Directional $\mathcal{H}^{2}$-Matrices for Helmholtz Problems
with Complex Frequency}
\author{Steffen B\"{o}rm\thanks{Institut f\"{u}r Informatik,
Christian-Albrechts-Universit\"{a}t zu Kiel, 24118 Kiel, Germany, e-mail:
\texttt{sb\texttt{@}informatik.uni-kiel.de}}
\and M. Lopez-Fernandez\thanks{Dipartamento di Matematica, Sapienza University of
Rome, Piazzale Aldo Moro 5, 00185 Roma, Italy, e-mail:
\texttt{lopez@mat.uniroma1.it}}
\and S.A. Sauter\thanks{Institut f\"{u}r Mathematik, Universit\"{a}t Z\"{u}rich,
Winterthurerstr 190, CH-8057 Z\"{u}rich, Switzerland, e-mail:
\texttt{stas@math.uzh.ch}}}
\maketitle

\begin{abstract}
The sparse approximation of high-frequency Helmholtz-type integral operators
has many important physical applications such as problems in wave propagation
and wave scattering. The discrete system matrices are huge and densely
populated; hence their sparse approximation is of outstanding importance. In
our paper we will generalize the directional $\mathcal{H}^{2}$-matrix
techniques from the \textquotedblleft pure\textquotedblright\ Helmholtz
operator $\mathcal{L}u=-\Delta u+\zeta^{2}u$ with $\zeta=-\operatorname*{i}k$,
$k\in\mathbb{R}$, to general complex frequencies $\zeta\in\mathbb{C}$ with
$\operatorname{Re}\zeta>0$. In this case, the fundamental solution decreases
exponentially for large arguments. We will develop a new admissibility
condition which contains $\operatorname{Re}\zeta$ in an explicit way and
introduce the approximation of the integral kernel function on admissible
blocks in terms of frequency-dependent \textit{directional expansion
functions}. We develop an error analysis which is explicit with respect to the
expansion order and with respect to $\operatorname{Re}\zeta$ and
$\operatorname{Im}\zeta$. This allows to choose the \textit{variable
}expansion order in a quasi-optimal way depending on $\operatorname{Re}\zeta$
but independent of, possibly large, $\operatorname{Im}\zeta$. The complexity
analysis is explicit with respect to $\operatorname{Re}\zeta$ and
$\operatorname{Im}\zeta$ and shows how higher values of $\operatorname{Re}%
\zeta$ reduce the complexity. In certain cases, it even turns out that the
discrete matrix can be replaced by its nearfield part.

Numerical experiments illustrate the sharpness of the derived estimates and
the efficiency of our sparse approximation.

\end{abstract}

\vspace{1em} \noindent\textbf{Keywords:} Helmholtz equation in lossy media,
hierarchical matrices, boundary integral operator

\vspace{1em} \noindent\textbf{Mathematics Subject Classification (2000):}
35J05, 65D05, 65N38, 41A10, 65N12

\section{Introduction}

The numerical simulation of many physical problems involves the solution of
large linear systems as a partial step in the overall algorithm. For
real-world applications the dimension of the system is huge, e.g., of order
$10^{6}-10^{10}$, which rules out exact elimination methods based, e.g., on
Gauss or Cholesky decompositions. Instead, iterative solvers are employed
which require a matrix-vector multiplication in each iteration step. If
non-local (integral) operators are involved the system matrices are fully
populated and a) the computation of the entries of the system matrix and b)
the matrix-vector multiplications typically are the bottlenecks in the
solution algorithms.

Since the mid 1980ies, the development of compression algorithms for\ densely
populated matrices related to the numerical discretization of non-local
operators has become an important topic in numerical analysis and scientific
computing. The fast multipole method has been developed in \cite{Rokhlin1} for
evaluating discrete Coulomb potentials. Panel-clustering methods have been
introduced first for collocation methods (see \cite{Nowak1}, \cite{NowakNum})
and were extended to Galerkin methods in \cite{SauterDiss}, \cite{Sautertsch}.
The idea of cluster methods have been generalized to a more algebraic setting
and led to the hierarchical $\mathcal{H}$-matrices (see \cite{HMatrixPart1},
\cite{Hackbusch_hier_matrices_engl}). A second hierarchy has been introduced
in \cite{Bungartz00}, \cite{BoermBook} and the resulting matrices are denoted
as $\mathcal{H}^{2}$ matrices. Most of these methods are restricted to
non-oscillatory elliptic problems. For highly oscillatory Helmholtz problems,
compression algorithms for the arising non-local integral operators have been
developed since the early 1990ies, among them are high-frequency, \textit{fast
multipole methods} \cite{RokhlinHelmholtz}, \cite{banjai_fast1},
\cite{banjai_fast2}, \cite{Kac-2014}, \textit{butterfly schemes}
\cite{Michielssen_butterfly}, \cite{Demanet_butterfly1},
\cite{boerm_butterfly}, and \textit{directional methods} \cite{Brandt_Helm},
\cite{Engquist_Ying}, \cite{Messner_Schanz_Darve}, \cite{Bebendorf_Kuske_Venn}%
, \cite{BoermMelenk}, \cite{Boerm2015}. For a comparison of these methods we
refer to \cite{BoermMelenk} and \cite{boerm_butterfly}.

The existing literature is mostly concerned with the \textquotedblleft
pure\textquotedblright\ Helmholtz problem, i.e., the operator $\mathcal{L}%
_{\zeta}u:=-\Delta u+\zeta^{2}u$ for purely imaginary frequency $\zeta
\in\operatorname*{i}\mathbb{R}$ (exceptions are the papers \cite{banjai_fast1}%
, \cite{banjai_fast2}, \cite{Kac-2014}). In our paper we consider more general
frequencies $\zeta\in\mathbb{C}$ with $\operatorname{Re}\zeta\geq0$ and recall
different important applications where such frequencies occur. The important
difference to the \textit{pure} Helmholtz problem is that the fundamental
solution exhibits an exponential decay for $\operatorname{Re}\zeta>0$ for
large arguments. We generalize the directional $\mathcal{H}^{2}$-matrix
approach in \cite{BoermMelenk} to complex frequencies and introduce new
admissibility conditions which contain the real part of the complex frequency
in an explicit way. We introduce the directional approximation of the integral
kernel function on admissible blocks and derive estimates for the
approximation error which allow to select the control parameters in a
quasi-optimal way. It turns out that a \textit{variable} expansion order
(depending on $\operatorname{Re}\zeta$) for different blocks is advantageous
compared to a fixed approximation order. In fact, it turns out that for
$\operatorname{Re}\zeta\gtrsim h_{\mathcal{G}}^{-1}$ (where $h_{\mathcal{G}}$
is a characteristic mesh width of the underlying boundary element mesh
$\mathcal{G}$) the discrete boundary element matrices can be replaced by its
nearfield part.

\begin{figure}[tbp] \centering

\includegraphics[
height=5.0548in,
width=5.0548in
]%
{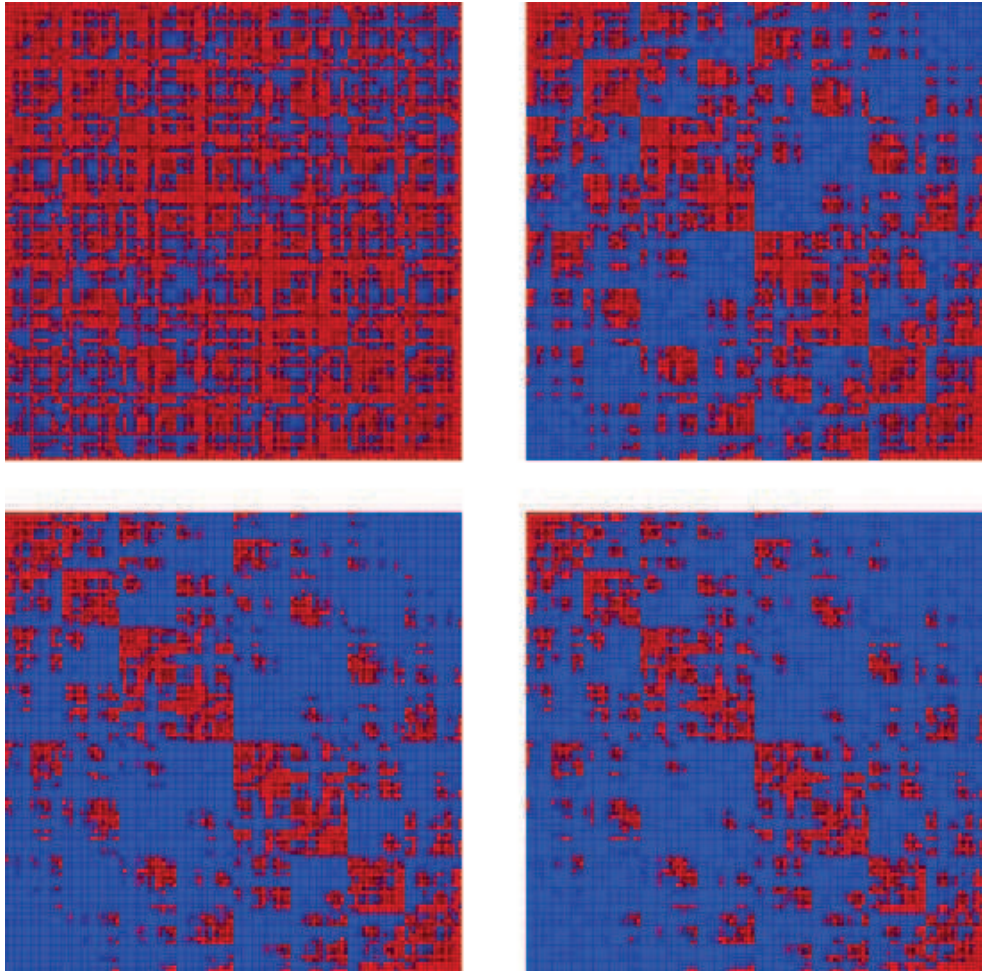}%

%\ $\caption{Sparsity pattern of the BEM matrices for the boundary integral
\caption{Sparsity pattern of the BEM matrices for the boundary integral
operator of the acoustic single layer potential. Non-admissible blocks are marked
by red and admissible blocks by blue. The ratio Re$\zeta$/Im$\zeta$ increases from
left to right and from top to bottom. On the left top, we consider the pure
Helmholtz problem, i.e., Im$\zeta$=0 while the ratios for the others are given by
1, 2, 3.\label{FigMat}}%
%TCIMACRO{\TeXButton{E}{\end{figure}}}%
%BeginExpansion
\end{figure}

The error estimates on the admissible blocks hold \textit{uniformly} with
respect to high oscillations. This has impact to the complexity analysis --
the compression rates benefit from a) the admissibility conditions which are
explicit in $\operatorname{Re}\zeta$, b) from the variable-order expansion,
and mostly from c) the frequency-explicit error estimates which allow to set
substantial parts of the system matrix to zero for large enough
$\operatorname{Re}\zeta$. In our numerical experiments (\S \ref{SecNumExp}),
we have applied our new admissibility condition and compression method to the
BEM matrix for the acoustic single layer potential. As an illustration we
depict in Figure \ref{FigMat} the dependence of the sparsity pattern on the
real and imaginary part of the wave number $\zeta$. As predicted by our
analysis the compression becomes stronger if the ratio $\operatorname{Re}%
\zeta/\operatorname{Im}\zeta$ increases.%
%TCIMACRO{\TeXButton{B}{\begin{figure}[tbp] \centering}}%
%BeginExpansion
%EndExpansion

The paper is structured as follows. In the next section, we will describe
three kinds of application where non-local Helmholtz-type integral operators
arise for general complex frequencies. In Section \ref{SecMethod} we formulate
the directional $\mathcal{H}^{2}$-matrix method with variable rank for general
complex frequencies. First, we introduce our new admissibility conditions and
then formulate the method in an algorithmic way. In Section \ref{SecAna}, we
estimate the error for the original integral kernel function being replaced by
the directional, variable order expansion on admissible matrix blocks. The
admissibility conditions along the error estimates form the basis for the
complexity analysis which is presented in Section \ref{SecComplexity}.
Finally, in Section \ref{SecNumExp} we report on the results of numerical
experiments which demonstrate the sharpness of our estimates.

\section{Setting\label{SecSet}}

In this section we will introduce three types of applications which lead to
Helmholtz-type equations at complex frequencies $\zeta\in\mathbb{C}$,
$\operatorname{Re}\zeta\geq0$.

\subsection{Helmholtz Equation with Decay}

Time harmonic wave propagation with decay arises in many applications such as,
e.g., in viscoelastodynamics for materials with damping (see, e.g.,
\cite{Achenbach}), in electromagnetics for wave propagation in lossy media
(see, e.g., \cite{Jackson02engl}), and in non-linear optics (see, e.g.,
\cite{sauter1996nonlinear}). In the simplest case such problems are modelled
by a Helmholtz equation with complex wave number.

\subsubsection{Variational Formulation}

Let $\Omega^{-}\subset\mathbb{R}^{3}$ be a bounded Lipschitz domain with
boundary $\Gamma$ and $\Omega^{+}:=\mathbb{R}^{3}\backslash\overline
{\Omega^{-}}$ its unbounded complement. For $\sigma\in\mathbb{R}$, let%
\[
\mathbb{C}_{\geq\sigma}:=\left\{  \zeta\in\mathbb{C}\mid\operatorname{Re}%
\zeta\geq\sigma\right\}  .
\]
Let the bilinear form $\left\langle \cdot,\cdot\right\rangle :\mathbb{C}%
^{3}\times\mathbb{C}^{3}\rightarrow\mathbb{C}$ be defined by $\left\langle
x,y\right\rangle =\sum_{j=1}^{3}x_{j}y_{j}$ so that the Euclidean norm is
given by $\left\Vert z\right\Vert =\left\langle z,\overline{z}\right\rangle
^{1/2}$. For a complex frequency $\zeta\in\mathbb{C}_{\geq0}$ and $\Omega
\in\left\{  \Omega^{-},\Omega^{+}\right\}  $ we consider the Helmholtz
equation subject to Dirichlet boundary conditions%
\begin{equation}%
\begin{array}
[c]{cc}%
-\Delta u+\zeta^{2}u=0 & \text{in }\Omega,\\
u=g_{\operatorname*{D}} & \text{on }\Gamma.
\end{array}
\label{condHelmEq}%
\end{equation}
If $\Omega=\Omega^{+}$ we also impose decay conditions at infinity%
\begin{equation}
\left\vert \dfrac{\partial w}{\partial r}+\zeta w\right\vert =Cr^{-1}%
\qquad\text{for }r:=\left\Vert x\right\Vert \rightarrow\infty.
\label{Sommerfeld3}%
\end{equation}
Here, $\partial_{r}$ denote the derivative in radial direction. We introduce
the acoustic Newton potential $G\left(  \zeta,\cdot\right)  $ by%
\begin{equation}
G\left(  \zeta,z\right)  :=\frac{\operatorname*{e}^{-\zeta\left\Vert
z\right\Vert }}{4\pi\left\Vert z\right\Vert }. \label{defkernels}%
\end{equation}
For the solution of (\ref{condHelmEq}) we employ an ansatz as an
\textit{acoustic single layer potential}%
\begin{equation}
\left(  S\left(  \zeta\right)  \varphi\right)  \left(  x\right)
:=\int_{\Gamma}G\left(  \zeta,y-x\right)  \varphi\left(  y\right)  d\Gamma
_{y}\qquad\forall x\in\Omega. \label{ansatz_slp}%
\end{equation}
To determine the \textit{unknown boundary density} $\varphi:\Gamma
\rightarrow\mathbb{C}$ we employ the Dirichlet boundary condition and the
continuity of the single layer operator up to the boundary. Let
\[
\left(  V\left(  \zeta\right)  \varphi\right)  \left(  x\right)
:=\int_{\Gamma}G\left(  \zeta,y-x\right)  \varphi\left(  y\right)  d\Gamma
_{y}\qquad\forall x\in\Gamma.
\]
Then, the strong formulation for the unknown density $\varphi$ is given by%
\begin{equation}
V\left(  \zeta\right)  \varphi=g_{\operatorname*{D}}\quad\text{on }\Gamma.
\label{strongform}%
\end{equation}
For the analysis of the \textit{boundary integral equation }and its Galerkin
discretization it is convenient to introduce the variational formulation. The
Sobolev spaces $H^{s}(\Gamma)$, $s\geq0,$ are defined in the usual way (see,
e.g., \cite{Hackell} or \cite{McLean00}) and the spaces with negative order
$s<0$ by duality. The norm is denoted by $\left\Vert \cdot\right\Vert
_{H^{s}\left(  \Gamma\right)  }$. The variational formulation of
(\ref{strongform}) is as follows: For given $g_{\operatorname*{D}}\in
H^{1/2}\left(  \Gamma\right)  $ find $\varphi\in H^{-1/2}\left(
\Gamma\right)  $ such that%
\begin{equation}
a_{\zeta}\left(  \varphi,\psi\right)  :=\left(  V\left(  \zeta\right)
\varphi,\psi\right)  =\left(  g_{\operatorname*{D}},\psi\right)  \qquad
\forall\psi\in H^{-1/2}\left(  \Gamma\right)  . \label{weakform}%
\end{equation}
Here $\left(  \cdot,\cdot\right)  $ denotes the continuous extension of the
$L^{2}\left(  \Gamma\right)  $ scalar product (with complex conjugation on the
second argument) to the anti-dual pairing on $H^{1/2}\left(  \Gamma\right)
\times H^{-1/2}\left(  \Gamma\right)  $, i.e.,%
\[
\left(  V\left(  \zeta\right)  \varphi,\psi\right)  =\int_{\Gamma}\int
_{\Gamma}G\left(  \zeta,y-x\right)  \varphi\left(  y\right)  \overline
{\psi\left(  x\right)  }d\Gamma_{y}d\Gamma_{x}.
\]

\begin{remark}
\label{RemExUniqu}Existence and uniqueness results for the solution of the
continuous problem (\ref{weakform}) \emph{for the case}\textit{ }%
$\operatorname{Re}\zeta>0$ are proved in \cite{Bam1}. For $\zeta
\in\operatorname*{i}\mathbb{R}$ it is well known that the operator $V\left(
\zeta\right)  $ is not invertible for discrete \textit{spurious frequencies}.
In this case, stabilized formulations (Brakhage-Werner or those proposed in
\cite{SauterBuffa}) cure this problem. We emphasize that $V\left(
\zeta\right)  $ appears in the stabilized formulations and its sparse
representation is still required.
\end{remark}

\subsubsection{Galerkin Discretization}

We consider the discretization of (\ref{weakform}) by a Galerkin boundary
element method. For a systematic introduction of boundary element methods we
refer, e.g., to the monograph \cite[Chap. 4]{SauterSchwab2010}. Let
$\mathcal{G=}\left\{  \tau_{i}:1\leq i\leq\tilde{M}\right\}  $ denote a
surface mesh of $\partial\Omega$, consisting of affine or possibly curved
triangles (called \textit{panels} in this context). As a convention the
triangles are (relatively) closed sets. For simplicity we assume that the
boundary element mesh $\mathcal{G}$ does not contain \textit{hanging nodes},
more precisely, that two non-identical triangles $\tau,t\in\mathcal{G}$ either
have a positive distance or their intersection is either a common edge or a
common vertex. For any $\tau\in\mathcal{G}$, there is a bijective
\textit{element map }$\chi_{\tau}:\widehat{\tau}\rightarrow\tau$ which maps
the reference element $\widehat{\tau}:=\operatorname*{conv}\left\{  \binom
{0}{0},\binom{1}{0},\binom{1}{1}\right\}  $ to the surface panel $\tau$; we
assume that this mapping is affine if $\tau$ is a plane triangle with straight
edges. In any case we assume that a common side $E$ of two adjacent triangles
$\tau,t\in\mathcal{G}$, are parametrized by $\chi_{\tau}$, $\chi_{t}$
\textquotedblleft in the same way\textquotedblright, i.e., $\chi_{\tau}%
^{-1}\left(  x\right)  =\gamma\circ\chi_{t}^{-1}\left(  x\right)  $ for all
$x\in E$ and a suitable affine mapping $\gamma:\widehat{\tau}\rightarrow
\widehat{\tau}$.

The finite-dimensional boundary element space of polynomial degree
$p\in\mathbb{N}_{0}$ and smoothness degree $m\in\left\{  -1,0\right\}  $
sub-ordinate to $\mathcal{G}$ is given by%
\[
S_{\mathcal{G}}^{p,m}:=\left\{  u\in L^{1}\left(  \Gamma\right)  \mid
\forall\tau\in\mathcal{G}\quad\left.  u\right\vert _{\overset{\circ}{\tau}%
}\circ\chi_{\tau}\in\mathbb{P}_{p}\right\}  \cap H^{m+1}\left(  \Gamma\right)
,
\]
where $\overset{\circ}{\tau}$ denotes the interior of $\tau$. If no confusion
is possible, we suppress the indices $p,m,\mathcal{G}$ and write $S$ short for
$S_{\mathcal{G}}^{p,m}$. The standard Lagrange nodal basis is denoted by
$b_{i}$, $i\in\mathcal{I}:=\left\{  1,\ldots,n\right\}  $, and depends as well
on $p$, $m$, $\mathcal{G}$.\ Finally we have%
\begin{equation}
S_{\mathcal{G}}^{p,m}=\operatorname{span}\left\{  b_{i}:1\leq i\leq n\right\}
\subset H^{-1/2}\left(  \partial\Omega\right)  . \label{DefVBEM}%
\end{equation}
The maximal mesh width is denoted by%
\[
h_{\mathcal{G}}:=\max\left\{  h_{\tau}:\tau\in\mathcal{G}\right\}
\quad\text{with\quad}h_{\tau}:=\operatorname{diam}\tau.
\]
For the mesh we define the \textit{shape-regularity }constants
$c_{\operatorname*{sr}}$ and $C_{\operatorname*{sr}}$ by\footnote{For a
measurable subset $\omega\subset\Gamma$ we denote by $\left\vert
\omega\right\vert $ the area measure of $\omega$.}%
\begin{equation}
c_{\operatorname*{sr}}:=\min_{\tau\in\mathcal{G}}\frac{\left\vert
\tau\right\vert }{h_{\tau}^{2}}\quad\text{and\quad}C_{\operatorname*{sr}%
}:=\max_{\tau\in\mathcal{G}}\frac{\left\vert \tau\right\vert }{h_{\tau}^{2}}.
\label{sr}%
\end{equation}
Another mesh parameter is%
\begin{equation}
h_{\min}:=\min\left\{  \operatorname*{dist}\left(  \tau,\tau^{\prime}\right)
:\forall\tau\in\mathcal{G\quad\forall\tau}^{\prime}\in\mathcal{G}\text{ with
}\tau\cap\tau^{\prime}=\emptyset\right\}  . \label{hmin}%
\end{equation}
We say that a boundary element mesh is \textit{quasi-uniform} if there exists
a constant $0<C_{\operatorname*{qu}}=\mathcal{O}\left(  1\right)  $ such that%
\begin{equation}
h_{\mathcal{G}}\leq C_{\operatorname*{qu}}h_{\min}. \label{quconst}%
\end{equation}
The number of panels is of the same order as the dimension of $S_{\mathcal{G}%
}^{p,m}$: there exist constants $C_{\operatorname*{loc},p}$, $C_{p}$ only
depending on the local polynomial degree $p$ such that for $m=-1,0$%
\begin{equation}
\max_{i\in\mathcal{I}}\sharp\left\{  j\in\mathcal{I}:\left\vert
\operatorname*{supp}b_{j}\cap\operatorname*{supp}b_{i}\right\vert \right\}
\leq C_{\operatorname*{loc},p}\quad\text{and\quad}C_{p}\sharp\mathcal{G\leq
}\dim S_{\mathcal{G}}^{p,m}=n. \label{cpconst}%
\end{equation}
The Galerkin discretization of equation (\ref{weakform}) is given by seeking
functions $\varphi_{S}\in S$ such that%
\begin{equation}
a_{\zeta}\left(  \varphi_{S},\psi\right)  =\left(  g_{\operatorname*{D}}%
,\psi\right)  \qquad\forall\psi\in S \label{weakform_discrete}%
\end{equation}
with $a_{\zeta}\left(  \cdot,\cdot\right)  $ as in (\ref{weakform}). By using
the basis $b_{i}$ we obtain a representation of this equation as a system of
linear equations. Let $\mathbf{K}\left(  \zeta\right)  =\left(  K_{i,j}\left(
\zeta\right)  \right)  _{i,j=1}^{n}\in\mathbb{C}^{n\times n}$ and
$\mathbf{r}=\left(  r_{i}\right)  _{i=1}^{n}$ be defined by%
\begin{equation}
K_{i,j}\left(  \zeta\right)  :=\left(  V\left(  \zeta\right)  b_{j}%
,b_{i}\right)  \quad\text{and\quad}r_{i}:=\left(  g_{\operatorname*{D}}%
,b_{i}\right)  ,\quad1\leq i,j\leq n. \label{defK}%
\end{equation}
The solution $%
%TCIMACRO{\TeXButton{boldphi}{\mbox{\boldmath$ \phi$}}}%
%BeginExpansion
\mbox{\boldmath$ \phi$}%
%EndExpansion
=\left(  \varphi_{i}\right)  _{i=1}^{n}$ of the linear system%
\[
\mathbf{K}\left(  \zeta\right)
%TCIMACRO{\TeXButton{boldphi}{\mbox{\boldmath$ \phi$}}}%
%BeginExpansion
\mbox{\boldmath$ \phi$}%
%EndExpansion
=\mathbf{r}%
\]
is then equivalent to the solution of (\ref{weakform_discrete}) via
\begin{equation}
\varphi_{S}=\sum_{i=1}^{n}\varphi_{i}b_{i}. \label{basisrep}%
\end{equation}

\begin{remark}
\quad

\begin{enumerate}
\item From \cite{Bam1} it follows that (\ref{weakform_discrete}) is well posed
for $\operatorname{Re}\zeta>0$ and, as a consequence, the matrix
$\mathbf{K}\left(  \zeta\right)  $ is invertible (cf. Rem. \ref{RemExUniqu}).

\item The matrix $\mathbf{K}\left(  \zeta\right)  $ is fully populated
containing, in general, $n^{2}$ non-zero entries. This is a major bottleneck
in a numerical realization of the boundary element method.

\item C\'{e}a's lemma can be applied to derive error estimates for the
solution $\varphi_{S}$. By using the results in \cite[Proof of Prop.
16]{LaSa2009} we obtain the quasi-optimal, frequency-explicit error estimate%
\[
\left\Vert \varphi-\varphi_{S}\right\Vert _{H^{-1/2}\left(  \Gamma\right)
}\leq C\frac{\left\vert \zeta\right\vert ^{3}}{\left(  \operatorname{Re}%
\zeta\right)  ^{2}}\inf_{\psi\in S}\left\Vert \varphi-\psi\right\Vert
_{H^{-1/2}\left(  \Gamma\right)  }.
\]

\item The theory deteriorates as $\operatorname{Re}\zeta\rightarrow0$ and this
is not an artifact. It is well known that the operator $V\left(  \zeta\right)
$ for purely imaginary wave number $\zeta$ is not injective for certain values
of $\zeta\in\operatorname*{i}\mathbb{R}$. Instead of (\ref{ansatz_slp}) one
often employs a combined double layer/ single layer ansatz, see, e.g.,
\cite{MelenkStab}, \cite{SauterBuffa}, and the resulting boundary integral
equation becomes well posed for all purely imaginary frequencies.
\end{enumerate}
\end{remark}

\subsection{Convolution Quadrature}

The linear homogeneous space-time wave equation can be transformed to
space-time boundary integral equations with retarded potentials. A popular
method for solving these equations is the convolution quadrature (CQ)
introduced in \cite{lub1}, \cite{lub2}. To circumvent the condition for the CQ
that the time steps must be constant, the $\operatorname*{gCQ}$ method has
been introduced in \cite{LoSau_1}, \cite{LoSau_alg}, \cite{LoSau_RK} to allow
for variable time steps. This method involves the numerical approximation of a
contour integral of the form%
\[
\frac{\Delta_{j}}{2\pi\operatorname*{i}}\int_{\mathcal{C}}\dfrac{V\left(
\zeta\right)  }{\prod_{\ell=j+1}^{N}(1-\Delta_{\ell}\zeta)}d\zeta,
\]
where $\mathcal{C}$ is a circle in the complex plane with midpoint $R>0$ and
radius $R$. The $j$-th time step is denoted by $\Delta_{j}$. For the numerical
evaluation of this contour integral, a quadrature method has been proposed in
\cite{LoSau_dd} which is of the form%
\[
\frac{\Delta_{j}}{2\pi\operatorname*{i}}{\sum_{\ell=1}^{N_{Q}}}w_{\ell}%
\dfrac{V\left(  \zeta_{\ell}\right)  }{\prod_{\ell=j+1}^{N}(1-\Delta_{\ell
}z_{\ell})}.
\]
Since the radius $R$ is large (proportional to the reciprocal minimal time
step) also the number of quadrature points is large and the numerical
realization requires the boundary element discretization of the operator
$V\left(  \zeta_{\ell}\right)  $ at all quadrature points on the contour
$\mathcal{C}$. Hence, also for this application one needs a sparse
approximation of the boundary element matrices $K\left(  \zeta_{\ell}\right)
$ for complex frequencies. We emphasize that for this application also the
case of non-resolved frequencies arises at certain quadrature points, i.e.,
the standard resolution condition $h_{\mathcal{G}}\left\vert
\operatorname*{Im}\zeta\right\vert \lesssim1$ is violated. The analysis of our
sparse approximation nicely reflects this fact: in certain cases (related to
the magnitude of $\operatorname{Re}\zeta$) only the nearfield part of the
system matrix has to be generated but not underresolved oscillations in the
farfield (cf. Remark \ref{RemTridiagCase}).

The need for the evaluation of such contour integrals also appears for the
original convolution quadrature method with constant time stepping since
non-local Helmholtz-type boundary element matrices have to be assembled in
many contour quadrature points. For the CQ method the resulting system matrix
is of block Toeplitz form and FFT-type techniques can be employed to reduce
the complexity with respect to the number of time points $N$ from
$\mathcal{O}\left(  N^{2}\right)  $ to almost $\mathcal{O}\left(  N\right)  $
(up to logarithmic terms). The combination of the FFT techniques in time and
sparse matrix techniques in space is far from trivial. In \cite{banjai_fast1},
\cite{banjai_fast2} such a fast multipole algorithm is introduced and
numerical experiments demonstrate the almost linear complexity (up to
logarithmic terms) with respect to the total number $Nn$ of unknowns; the
generalization to general complex frequencies of this fast multipole method is
presented and analyzed in \cite{Kac-2014}. The spatial compression algorithm
for the Helmholtz-type boundary element matrices is based on the
high-frequency multipole method which goes back to \cite{RokhlinHelmholtz} and
is different from the directional $\mathcal{H}^{2}$ matrices (for a comparison
of these compression methods we refer to \cite[\S 1]{BoermMelenk}). We expect
that our directional $\mathcal{H}^{2}$-matrix compression algorithm has the
potential to be used within the fast method described in \cite{banjai_fast1},
\cite{banjai_fast2} and has the advantage that a fully developed accuracy
analysis is available to select the control parameters in a quasi-optimal way.

\subsection{Limiting Absorbing Principle}

The transformation of the time-space linear wave equation $\left(
\partial_{t}^{2}-\Delta\right)  u=f$ to the frequency domain by a time
periodic ansatz leads to the Helmholtz equation of the form (\ref{condHelmEq})
with purely imaginary frequency $\zeta\in\operatorname*{i}\mathbb{R}$. This
equation is not solvable if \textquotedblleft$-\zeta^{2}$\textquotedblright%
\ is an eigenvalue of the (negative) Laplacian with Dirichlet boundary
conditions. For theoretical as well as for practical reasons (see, e.g.,
\cite{limit_absorb_princ_nguyen}, \cite{Vuik_Helm_Solv}) it can be useful to
\textquotedblleft add some absorption\textquotedblright\ to this equation and
to consider the equation%
\[%
\begin{array}
[c]{cc}%
-\Delta u+\left(  \zeta^{2}-\operatorname*{i}\varepsilon\right)  u=0 &
\text{in }\Omega,\\
u=g_{\operatorname*{D}} & \text{on }\Gamma
\end{array}
\]
for a small positive parameter $\varepsilon$. This equation is solvable for
all frequencies $\zeta\in\operatorname*{i}\mathbb{R}$ and one can employ the
Galerkin boundary element for its discretization. The discretization matrix is
given by $\mathbf{K}\left(  \tilde{\zeta}\right)  $ (cf. (\ref{defK})) for the
choice $\tilde{\zeta}=\sqrt{\zeta^{2}-\operatorname*{i}\varepsilon
}=a+\operatorname*{i}b$ with $a:=\sqrt{\frac{\zeta^{2}+\sqrt{\zeta
^{4}+\varepsilon^{2}}}{2}}>0$ and $b:=-\varepsilon/\left(  2a\right)  $. This
implies$\ \tilde{\zeta}\in\mathbb{C}_{>0}$. Hence, also in this case sparse
matrix techniques which are applicable to $\mathbf{K}\left(  \zeta\right)  $,
$\zeta\in\mathbb{C}_{>0}$ are important for the Helmholtz equation with
artificially added absorption.

\section{Directional $\mathcal{H}^{2}$ Matrices for Helmholtz Equations with
Decay\label{SecMethod}}

Directional $\mathcal{H}^{2}$ matrices have been introduced in
\cite{Brandt_Helm}, \cite{Engquist_Ying}, \cite{Messner_Schanz_Darve},
\cite{Bebendorf_Kuske_Venn}, \cite{BoermMelenk}, \cite{Boerm2015} for the
high-frequency Helmholtz problems for purely imaginary frequency.

\subsection{Directionally Admissible Partitionings}

We generalize this method to general complex frequencies and analyze its
accuracy and complexity explicitly with respect to the real and imaginary part
of the wave number. To formulate the algorithm we first introduce some notation.

As a basis for the boundary element space we have chosen Lagrange basis
functions $b_{i}$ which have local support $\omega_{i}:=\operatorname*{supp}%
b_{i}$. We collect the set of indices $1\leq i\leq n$ in the set $\mathcal{I}$
so that $\sharp\mathcal{I}=n$.

\begin{definition}
\label{DefCluster}For a given set of degrees of freedom $\mathcal{I}$,
the\emph{ cluster tree }$\mathcal{T}_{\mathcal{I}}$ is a \emph{labeled tree}
which satisfies:

\begin{enumerate}
\item the label $\hat{t}$ of each node $t\in\mathcal{T}_{\mathcal{I}}$ is a
subset of the index set $\mathcal{I}$,

\item the \emph{root} $r\in\mathcal{T}_{\mathcal{I}}$ of the tree is assigned
$\hat{r}=\mathcal{I}$;

\item for all $t\in\mathcal{T}_{\mathcal{I}}$ there exists

\begin{enumerate}
\item either a set of nodes $\operatorname*{sons}\left(  t\right)  $ denoted
by \emph{sons of }$t$ which satisfies: $\hat{t}=%
%TCIMACRO{\dbigcup \limits_{s\in\operatorname*{sons}\left(  t\right)  }}%
%BeginExpansion
{\displaystyle\bigcup\limits_{s\in\operatorname*{sons}\left(  t\right)  }}
%EndExpansion
\hat{s}$ and, for all $t_{1}$, $t_{2}\in\operatorname*{sons}\left(  t\right)
$, it holds either $t_{1}=t_{2}$ or $\hat{t}_{1}\cap\hat{t}_{2}=\emptyset$.

\item or $t$ is called a \emph{leaf}. The set of leaf clusters is%
\[
\mathcal{L}_{\mathcal{I}}:=\left\{  t\in\mathcal{T}_{\mathcal{I}}:t\text{ is a
leaf}\right\}  .
\]
Vice versa $t=\operatorname*{father}\left(  t^{\prime}\right)  $ is the
\emph{father} of $t^{\prime}\in\operatorname*{sons}\left(  t\right)  $;
\end{enumerate}

\item[4.] With each cluster, an axis-parallel \emph{bounding box} $B_{t}$ is
associated which satisfies
\begin{equation}
\omega_{t}:=%
%TCIMACRO{\dbigcup \limits_{i\in\hat{t}}}%
%BeginExpansion
{\displaystyle\bigcup\limits_{i\in\hat{t}}}
%EndExpansion
\omega_{i}\subset B_{t}. \label{clusterbox}%
\end{equation}
The \textit{center }$M_{t}$ of a cluster $t$ is defined as the barycenter of
$B_{t.}$
\end{enumerate}

The \emph{level} of a cluster is given recursively by $\operatorname*{level}%
\left(  r\right)  :=0$ and $\operatorname*{level}\left(  t^{\prime}\right)
:=\operatorname*{level}\left(  t\right)  +1$ for all $t\in\mathcal{T}%
_{\mathcal{I}}\backslash\mathcal{L}_{\mathcal{I}}$ and $t^{\prime}%
\in\operatorname*{sons}\left(  t\right)  $. The \emph{depth} of a cluster tree
is $\operatorname*{depth}\left(  \mathcal{T}_{\mathcal{I}}\right)
:=\max\left\{  \operatorname*{level}\left(  t\right)  :t\in\mathcal{T}%
_{\mathcal{I}}\right\}  $ and $\mathcal{T}_{\ell}:=\left\{  t\in
\mathcal{T}_{\mathcal{I}}\mid\operatorname*{level}\left(  t\right)
=\ell\right\}  $. The maximal \emph{cluster diameter} of level $\ell$ is%
\begin{equation}
\delta_{\ell}:=\max\left\{  \operatorname*{diam}B_{t}:t\in\mathcal{T}_{\ell
}\right\}  . \label{maxcldiam}%
\end{equation}

\end{definition}

A natural choice for the bounding box $B_{t}$ is the minimal box such that
(\ref{clusterbox}) holds but we do not restrict to this choice. However, we
assume that there exist positive constants $c_{\operatorname*{vol}}$,
$C_{\operatorname*{vol}}$ such that%
\begin{equation}
c_{\operatorname*{vol}}\left\vert \omega_{t}\right\vert \leq
\operatorname*{diam}\nolimits^{2}B_{t}\leq C_{\operatorname*{vol}}\left\vert
\omega_{t}\right\vert . \label{Cwidth}%
\end{equation}

Algorithms for building cluster trees from index sets $\mathcal{I}$
corresponding to boundary element basis functions can be found, e.g., in
\cite{Sauter00}, \cite{Hackbusch_hier_matrices_engl}.

The clusters allow via the geometric correspondence (\ref{clusterbox}) to
identify pairs of regions $B_{t},B_{s}\subset\Gamma$ -- associated to pairs of
clusters $\left(  t,s\right)  $ -- where the kernel function can be
approximated by a separable expansion:%
\[
G\left(  \zeta,y-x\right)  \approx\sum_{\nu=1}^{k}\sum_{\mu=1}^{k}\gamma
_{\nu,\mu}\Phi_{\nu}^{t}\left(  x\right)  \Psi_{\mu}^{s}\left(  y\right)
\qquad\forall\left(  x,y\right)  \in B_{t}\times B_{s},
\]
i.e., an expansion where $x$ and $y$ appear only in a factorized way. The
number $k$ is denoted as the \textit{rank} of the separable expansion. To
identify these regions we employ an admissibility condition which will be
introduced next. It will turn out from our analysis that for a pair $\left(
t,s\right)  $ of \textit{admissible} clusters the kernel function can be
approximated by a separable expansion.

\begin{definition}
[directional admissibility condition for complex frequencies]For $%
%TCIMACRO{\TeXButton{boldeta}{\mbox{\boldmath$ \eta$}}}%
%BeginExpansion
\mbox{\boldmath$ \eta$}%
%EndExpansion
=\left(  \eta_{i}\right)  _{i=1}^{3}\in\mathbb{R}_{>0}^{3}$ a pair of clusters
$t,s\in\mathcal{T}_{\mathcal{I}}$ and a direction $c\in\mathbb{S}_{2}$ are $%
%TCIMACRO{\TeXButton{boldeta}{\mbox{\boldmath$ \eta$}}}%
%BeginExpansion
\mbox{\boldmath$ \eta$}%
%EndExpansion
$\emph{-admissible} with respect to a complex frequency $\zeta\in
\mathbb{C}_{>0}$ if they satisfy the following three conditions:%
%TCIMACRO{\TeXButton{adm}{\begin{subequations}
%\label{adm}
%\end{subequations}}}%
%BeginExpansion
\begin{subequations}
\label{adm}
\end{subequations}%
%EndExpansion%
\begin{align}
\left\vert \operatorname{Im}\zeta\right\vert \left\Vert \frac{M_{t}-M_{s}%
}{\left\Vert M_{t}-M_{s}\right\Vert }-c\right\Vert  &  \leq\frac{\eta_{1}%
}{\max\left\{  \operatorname*{diam}\left(  B_{t}\right)  ,\operatorname*{diam}%
\left(  B_{s}\right)  \right\}  },\tag{%
%TCIMACRO{\TeXButton{adm}{\ref{adm}}}%
%BeginExpansion
\ref{adm}%
%EndExpansion
a}\label{adm_angle}\\
\max\left\{  \operatorname*{diam}\left(  B_{t}\right)  ,\operatorname*{diam}%
\left(  B_{s}\right)  \right\}   &  \leq\eta_{2}\operatorname*{dist}\left(
B_{t},B_{s}\right)  ,\tag{%
%TCIMACRO{\TeXButton{adm}{\ref{adm}}}%
%BeginExpansion
\ref{adm}%
%EndExpansion
b}\label{adm_sing}\\
\left\vert \operatorname{Im}\zeta\right\vert \max\left\{  \operatorname*{diam}%
\nolimits^{2}\left(  B_{t}\right)  ,\operatorname*{diam}\nolimits^{2}\left(
B_{s}\right)  \right\}   &  \leq\max\left\{  \eta_{2},\eta_{3}\left(
\operatorname{Re}\zeta\right)  \operatorname*{dist}\left(  B_{t},B_{s}\right)
\right\}  \operatorname*{dist}\left(  B_{t},B_{s}\right)  . \tag{%
%TCIMACRO{\TeXButton{adm}{\ref{adm}}}%
%BeginExpansion
\ref{adm}%
%EndExpansion
c}\label{adm_paranew}%
\end{align}

\end{definition}

The algorithm for generating a minimal partition $\mathcal{P}$ of
$\mathcal{I}\times\mathcal{I}$ by admissible and non-admissible blocks is of
divide-and-conquer type.

\begin{remark}
\label{Remdeltamin}We will need the minimal distance between the clusters of
admissible blocks and set%
\[
\delta_{\min}:=\min\left\{  \operatorname*{dist}\left(  B_{t},B_{s}\right)
:\left(  t,s\right)  \text{ is admissible}\right\}  .
\]
We have by (\ref{adm_sing})%
\[
\delta_{\min}\geq\frac{1}{\eta_{2}}\min\left\{  \operatorname*{diam}\left(
B_{t}\right)  :t\in\mathcal{T}_{\mathcal{I}}\right\}  \geq\frac{h_{\min}}%
{\eta_{2}}%
\]
with $h_{\min}$ as in (\ref{hmin}).
\end{remark}

\begin{algorithm}
\label{Algodivide}The minimal, $%
%TCIMACRO{\TeXButton{boldeta}{\mbox{\boldmath$ \eta$}}}%
%BeginExpansion
\mbox{\boldmath$ \eta$}%
%EndExpansion
$-admissible block partitioning $\mathcal{P}$ of $\mathcal{I}\times
\mathcal{I}$ is obtained as the result of the recursive procedure
\textbf{divide}$\left(  \left(  r,r\right)  ,\emptyset\right)  $ defined
by\medskip\ (see \cite{NowakNum})
\end{algorithm}

\textbf{procedure} \textbf{divide}$\left(  b,\mathcal{P}\right)  $;

\textbf{begin \{}notation: $b=\left(  t,s\right)  $ for $t,s\in\mathcal{T}%
_{\mathcal{I}}$\}

\quad\textbf{if} $\operatorname*{sons}\left(  t\right)  =\emptyset$\textbf{
or} $\operatorname*{sons}\left(  s\right)  =\emptyset$ \textbf{then begin}

$\quad\quad\operatorname*{sons}\left(  b\right)  \leftarrow\emptyset$;

$\quad\quad\mathcal{P}\leftarrow\mathcal{P}\cup\left\{  b\right\}  $

$\quad$\textbf{end else if} ($b$ is admissible) \textbf{then begin}

$\quad\quad\operatorname*{sons}\left(  b\right)  \leftarrow\emptyset$:

$\quad\quad\mathcal{P}\leftarrow\mathcal{P}\cup\left\{  b\right\}  $

$\quad$\textbf{end else} \textbf{for} $t^{\prime}\in\operatorname*{sons}%
\left(  t\right)  $, $s^{\prime}\in\operatorname*{sons}\left(  s\right)  $
\textbf{do begin}

$\quad\quad b^{\prime}\leftarrow\left(  t^{\prime},s^{\prime}\right)  $;

$\quad\quad$\textbf{divide}($b^{\prime},\mathcal{P}$)

\textbf{end end\medskip}

We split the covering $\mathcal{P}=\mathcal{P}_{\operatorname*{near}}%
\cup\mathcal{P}_{\operatorname*{far}}$ with%
\begin{equation}
\mathcal{P}_{\operatorname*{far}}:=\left\{  b\in\mathcal{P}\mid b\text{ is
admissible}\right\}  \quad\text{and\quad}\mathcal{P}_{\operatorname*{near}%
}:=\mathcal{P}\backslash\mathcal{P}_{\operatorname*{far}}.\label{defPnearPfar}%
\end{equation}
For a cluster $t\in\mathcal{T}_{\mathcal{I}}$, we define the set of
\textit{left} and \textit{right} \textit{partners} by%
\begin{equation}%
\begin{array}
[c]{ll}%
\mathcal{P}_{\operatorname*{left}}^{\operatorname*{near}}\left(  t\right)
:=\left\{  s:\left(  s,t\right)  \in\mathcal{P}_{\operatorname*{near}%
}\right\}  , & \mathcal{P}_{\operatorname*{right}}^{\operatorname*{near}%
}\left(  t\right)  :=\left\{  s:\left(  t,s\right)  \in\mathcal{P}%
_{\operatorname*{near}}\right\}  ,\\
\mathcal{P}_{\operatorname*{left}}^{\operatorname*{far}}\left(  t\right)
:=\left\{  s:\left(  s,t\right)  \in\mathcal{P}_{\operatorname*{far}}\right\}
, & \mathcal{P}_{\operatorname*{right}}^{\operatorname*{far}}\left(  t\right)
:=\left\{  s:\left(  t,s\right)  \in\mathcal{P}_{\operatorname*{far}}\right\}
,\\
\mathcal{P}_{\operatorname*{left}}\left(  t\right)  :=\mathcal{P}%
_{\operatorname*{left}}^{\operatorname*{near}}\left(  t\right)  \cup
\mathcal{P}_{\operatorname*{left}}^{\operatorname*{far}}\left(  t\right)  , &
\mathcal{P}_{\operatorname*{right}}\left(  t\right)  :=\mathcal{P}%
_{\operatorname*{right}}^{\operatorname*{near}}\left(  t\right)
\cup\mathcal{P}_{\operatorname*{right}}^{\operatorname*{far}}\left(  t\right)
\end{array}
\label{leftright}%
\end{equation}

\begin{remark}
\label{RemLevels}We have not assumed that $\mathcal{T}_{\mathcal{I}}$ is a
balanced tree\footnote{A balanced tree is a tree where $\mathcal{L}%
_{\mathcal{I}}=\mathcal{T}_{L}$ for $L:=\operatorname*{depth}\left(
\mathcal{T}_{\mathcal{I}}\right)  $.}. However, the definition of the sons of
a block $b=\left(  t,s\right)  $ imply that each block $b=\left(  t,s\right)
\in\mathcal{P}$ consists of clusters $t,s\in\mathcal{T}_{\mathcal{I}}$ which
have the same level in the cluster tree $\operatorname*{level}\left(
t\right)  =\operatorname*{level}\left(  s\right)  $ and we set
$\operatorname*{level}\left(  b\right)  :=\operatorname*{level}\left(
t\right)  $. As a consequence, we have $\operatorname*{depth}\mathcal{P}%
=\operatorname*{depth}\mathcal{T}_{\mathcal{I}}$.
\end{remark}

\subsection{Approximation of the Kernel Function}

Next we explain the approximation of the kernel function $G\left(  \zeta
,\cdot\right)  $ for the single layer boundary integral operator (cf.
(\ref{defkernels})) for complex frequencies. For some unit vector
$c\in\mathbb{S}_{2}$ we write\footnote{The bilinear form $\left\langle
\cdot,\cdot\right\rangle :\mathbb{C}^{r}\times\mathbb{C}^{r}\rightarrow
\mathbb{C}$ is defined by $\left\langle x,y\right\rangle =\sum_{j=1}^{r}%
x_{j}y_{j}$.}%
\begin{equation}
G\left(  \zeta,z\right)  :=\frac{\operatorname*{e}^{-\zeta\left\Vert
z\right\Vert }}{4\pi\left\Vert z\right\Vert }=\operatorname*{e}%
\nolimits^{-\operatorname*{i}\left(  \operatorname*{Im}\zeta\right)
\left\langle z,c\right\rangle }G_{c}\left(  \zeta,z\right)  \label{degGnew0}%
\end{equation}
with%
\begin{equation}
G_{c}\left(  \zeta,z\right)  :=\operatorname*{e}\nolimits^{-\left(
\operatorname{Re}\zeta\right)  \left\Vert z\right\Vert }\frac
{\operatorname*{e}\nolimits^{-\operatorname*{i}\left(  \operatorname*{Im}%
\zeta\right)  \left(  \left\Vert z\right\Vert -\left\langle z,c\right\rangle
\right)  }}{4\pi r}. \label{defcnew}%
\end{equation}
Let $b=\left(  t,s\right)  $ be an admissible block. We approximate this
kernel function on $B_{t}\times B_{s}$ by%
\[
\tilde{G}_{b}\left(  \zeta,\cdot\right)  :=\operatorname*{e}%
\nolimits^{-\operatorname*{i}\left(  \operatorname{Im}\zeta\right)
\left\langle \cdot,c\right\rangle }\mathfrak{I}_{b}\left(  G_{c}\left(
\zeta,\cdot\right)  \right)  ,
\]
where $\mathfrak{I}_{b}$ denotes the tensor \v{C}eby\v{s}ev interpolation on
$B_{t}\times B_{s}$ with polynomials of maximal degree $m$. The degree $m$ as
well as the direction $c$ depend on the block $b$, i.e., $m=m\left(  b\right)
$ and $c=c\left(  b\right)  $. The choice of $c\left(  b\right)  $ will be
explained next. From the error analysis/admissibility condition it follows
that an ideal choice is%
\[
c=\frac{M_{t}-M_{s}}{\left\Vert M_{t}-M_{s}\right\Vert }.
\]
However, for efficiency reasons we restrict the number of possible choices of
directions to a finite set $\mathcal{D}_{\ell}$ which will depend on the level
$\ell$ of a block. Recall the definition of the maximal cluster diameter
$\delta_{\ell}$ on level $\ell$ (cf. (\ref{maxcldiam})). The finite set
$\mathcal{D}_{\ell}\subset\mathbb{S}_{2}$ has to satisfy by (by
(\ref{adm_angle}))%
\[
\left\vert \operatorname{Im}\zeta\right\vert \sup_{e\in\mathbb{S}_{2}}%
\inf_{c\in\mathcal{D}_{\ell}}\left\Vert e-c\right\Vert \leq\frac{\eta_{1}%
}{\delta_{\ell}},
\]
which guarantees that, for any block $b$ which satisfies (\ref{adm_sing}) and
(\ref{adm_paranew}), there exists a direction $c\left(  b\right)
\in\mathcal{D}_{\operatorname*{level}\left(  b\right)  }$ such that
(\ref{adm_angle}) is also satisfied and the block is admissible. There are
various methods to construct such sets of directions. Here, we choose the
construction as explained in \cite[Rem. 3]{borm2018hybrid}. Since
$\delta_{\ell}$ is smaller for finer levels ($\ell$ large) we may conclude
that the cardinality of the set $\mathcal{D}_{\ell}$ increases for larger blocks.

Next we explain the choice of $m$. Let $b=\left(  t,s\right)  $ be an
admissible block. In Section \ref{SecAna} we will prove the estimate%
\[
\left\Vert G\left(  \zeta,\cdot\right)  -\tilde{G}_{b}\left(  \zeta
,\cdot\right)  \right\Vert _{\infty,B_{t}-B_{s}}\leq\frac{C_{0}%
\operatorname*{e}^{-\sigma\left(  \operatorname{Re}\zeta\right)
\operatorname*{dist}\left(  B_{t},B_{s}\right)  }}{4\pi\operatorname*{dist}%
\left(  B_{t},B_{s}\right)  }\rho_{0}^{-m}%
\]
for some $\sigma,C_{0}>0$ and $\rho_{0}>1$. Here $B_{t}-B_{s}:=\left\{
y-x\mid\left(  x,y\right)  \in B_{t}\times B_{s}\right\}  $. If we aim for a
constant error $\varepsilon$ on each block\footnote{To simplify the
calculations we restrict to $0<\varepsilon\leq\min\left\{  \operatorname*{e}%
^{-1},\frac{h_{\min}}{\eta_{2}}\right\}  $ so that $\log\frac{1}{\varepsilon
}\geq\log\frac{1}{\operatorname*{dist}\left(  B_{t},B_{s}\right)  }$ for all
admissible blocks $b=\left(  t,s\right)  $ (cf. Rem. \ref{Remdeltamin}).} the
condition
\[
\frac{C_{0}\operatorname*{e}^{-\sigma\left(  \operatorname{Re}\zeta\right)
\operatorname*{dist}\left(  B_{t},B_{s}\right)  }}{4\pi\operatorname*{dist}%
\left(  B_{t},B_{s}\right)  }\rho_{0}^{-m}\leq\varepsilon
\]
leads to a dependence of $m$ on $\operatorname{Re}\zeta$, on the block
$b=\left(  t,s\right)  $ and on $\varepsilon$ of the form%
\begin{equation}
\tilde{m}_{b}:=\left\{
\begin{array}
[c]{ll}%
\left\lceil c_{0}\log\frac{1}{\varepsilon}-\tilde{\sigma}\left(
\operatorname{Re}\zeta\right)  \operatorname*{dist}\left(  B_{t},B_{s}\right)
\right\rceil  & \text{if }c_{0}\log\frac{1}{\varepsilon}\geq\tilde{\sigma
}\left(  \operatorname{Re}\zeta\right)  \operatorname*{dist}\left(
B_{t},B_{s}\right)  ,\\
-1 & \text{otherwise}%
\end{array}
\right.  \label{defmbtilde}%
\end{equation}
for positive constants $c_{0}$, $\tilde{\sigma}>0$. As a convention,
$\tilde{m}_{b}=-1$ means that the kernel function is replaced by the zero
function on this block. To keep the algorithm simple we will restrict to a
single expansion order per level $\ell$ by setting%
\[
\tilde{m}_{\ell}:=\max\left\{  \tilde{m}_{b}:b\in\mathcal{P}%
_{\operatorname*{far}}\quad\text{with\quad}\operatorname*{level}\left(
b\right)  =\ell\right\}  .
\]
In order to get the second (functional) hierarchy (besides the geometric
cluster hierarchy) which allows to represent the expansion on larger clusters
by an expansion on smaller clusters it is necessary that the sequence of
expansion orders $\left(  m_{\ell}\right)  _{\ell=0}^{\operatorname*{depth}%
\left(  \mathcal{T}_{\mathcal{I}}\right)  }$ is increasing towards the leaves.
This is guaranteed by the recursive definition%
\begin{equation}
m\left(  b\right)  :=m_{\ell}:=\left\{
\begin{array}
[c]{ll}%
\tilde{m}_{0} & \ell=0,\\
\max\left\{  m_{\ell-1},\tilde{m}_{\ell}\right\}  & \ell\geq1
\end{array}
\right.  \qquad\forall b\in\mathcal{P}_{\operatorname*{far}}\text{ with
}\operatorname*{level}\left(  b\right)  =\ell. \label{defmb}%
\end{equation}

We pass the approximation orders from the admissible blocks onto the clusters
via (cf. Remark \ref{RemLevels})%
\[
m_{t}:=m_{\operatorname{level}\left(  t\right)  }.
\]
The (first) approximation of the kernel function on an admissible block
$b=\left(  t,s\right)  $ with expansion order $m=m\left(  b\right)  $ and
$c=c\left(  b\right)  $ can be written in the form%
\[
\tilde{G}_{b}\left(  \zeta,y-x\right)  :=\sum_{%
%TCIMACRO{\TeXButton{boldmu}{\mbox{\boldmath$ \mu$}}}%
%BeginExpansion
\mbox{\boldmath$ \mu$}%
%EndExpansion
,%
%TCIMACRO{\TeXButton{boldnu}{\mbox{\boldmath$ \nu$}}}%
%BeginExpansion
\mbox{\boldmath$ \nu$}%
%EndExpansion
\in\mathbb{N}_{t}}\gamma_{%
%TCIMACRO{\TeXButton{boldmu}{\mbox{\boldmath$ \mu$}}}%
%BeginExpansion
\mbox{\boldmath$ \mu$}%
%EndExpansion
,%
%TCIMACRO{\TeXButton{boldnu}{\mbox{\boldmath$ \nu$}}}%
%BeginExpansion
\mbox{\boldmath$ \nu$}%
%EndExpansion
,c}^{b}\left(  \zeta\right)  \tilde{\Phi}_{%
%TCIMACRO{\TeXButton{boldmu}{\mbox{\boldmath$ \mu$}}}%
%BeginExpansion
\mbox{\boldmath$ \mu$}%
%EndExpansion
,c}^{t}\left(  \zeta,x\right)  \overline{\tilde{\Phi}_{%
%TCIMACRO{\TeXButton{boldnu}{\mbox{\boldmath$ \nu$}}}%
%BeginExpansion
\mbox{\boldmath$ \nu$}%
%EndExpansion
,c}^{s}\left(  \zeta,y\right)  },
\]
where we employ the following notation: The index set $\mathbb{N}_{t}$ is
given by%
\begin{equation}
\mathbb{N}_{t}:=\left\{
%TCIMACRO{\TeXButton{boldmu}{\mbox{\boldmath$ \mu$}}}%
%BeginExpansion
\mbox{\boldmath$ \mu$}%
%EndExpansion
\in\mathbb{N}_{0}^{3}\mid0\leq\mu_{i}\leq m_{t},\quad1\leq i\leq3\right\}
\label{defnt}%
\end{equation}
and we denote by%
\[
k_{\operatorname{level}\left(  t\right)  }:=k_{t}:=\sharp\mathbb{N}_{t}%
\]
the \textit{rank }of the expansion for $t\in\mathcal{T}_{\operatorname{level}%
\left(  t\right)  }$. Note that $\mathbb{N}_{t}=\emptyset$ for $m_{t}=-1$ so
that $k_{t}=0$. Let $\hat{\xi}_{i,m}$, $0\leq i\leq m$, denote the
\v{C}eby\v{s}ev nodal points on the unit interval $\left(  -1,1\right)  $ and
let $\hat{L}_{i,m}$ be the corresponding Lagrange polynomials. The tensor
versions are given, for $%
%TCIMACRO{\TeXButton{boldmu}{\mbox{\boldmath$ \mu$}}}%
%BeginExpansion
\mbox{\boldmath$ \mu$}%
%EndExpansion
\in\mathbb{N}_{0}^{3}$, $0\leq\mu_{i}\leq m$, by $\hat{\xi}_{%
%TCIMACRO{\TeXButton{boldmu}{\mbox{\boldmath$ \mu$}}}%
%BeginExpansion
\mbox{\boldmath$ \mu$}%
%EndExpansion
,m}:=\left(  \hat{\xi}_{\mu_{1},m},\hat{\xi}_{\mu_{2},m},\hat{\xi}_{\mu_{3}%
,m}\right)  ^{\intercal}$ and $\hat{L}_{%
%TCIMACRO{\TeXButton{boldmu}{\mbox{\boldmath$ \mu$}}}%
%BeginExpansion
\mbox{\boldmath$ \mu$}%
%EndExpansion
,m}\left(  x\right)  =%
%TCIMACRO{\dprod \limits_{\ell=1}^{3}}%
%BeginExpansion
{\displaystyle\prod\limits_{\ell=1}^{3}}
%EndExpansion
\hat{L}_{\mu_{\ell},m}\left(  x_{\ell}\right)  $. For a box $B_{t}$, let
$\theta_{t}$ denote an affine pullback to the cube $\left(  -1,1\right)  ^{3}%
$. Then, the tensorized \v{C}eby\v{s}ev nodal points of order $m_{t}$ scaled
to the sides of $B_{t}$ are given by $\xi_{%
%TCIMACRO{\TeXButton{boldmu}{\mbox{\boldmath$ \mu$}}}%
%BeginExpansion
\mbox{\boldmath$ \mu$}%
%EndExpansion
,t}:=\theta_{t}\left(  \hat{\xi}_{%
%TCIMACRO{\TeXButton{boldmu}{\mbox{\boldmath$ \mu$}}}%
%BeginExpansion
\mbox{\boldmath$ \mu$}%
%EndExpansion
,m_{t}}\right)  $ and $L_{%
%TCIMACRO{\TeXButton{boldmu}{\mbox{\boldmath$ \mu$}}}%
%BeginExpansion
\mbox{\boldmath$ \mu$}%
%EndExpansion
}^{t}:=\hat{L}_{%
%TCIMACRO{\TeXButton{boldmu}{\mbox{\boldmath$ \mu$}}}%
%BeginExpansion
\mbox{\boldmath$ \mu$}%
%EndExpansion
,m_{t}}\circ\theta_{t}^{-1}$, for all $%
%TCIMACRO{\TeXButton{boldmu}{\mbox{\boldmath$ \mu$}}}%
%BeginExpansion
\mbox{\boldmath$ \mu$}%
%EndExpansion
\in\mathbb{N}_{t}$. The expansion functions are given by%
\[
\tilde{\Phi}_{%
%TCIMACRO{\TeXButton{boldmu}{\mbox{\boldmath$ \mu$}}}%
%BeginExpansion
\mbox{\boldmath$ \mu$}%
%EndExpansion
,c}^{t}\left(  \zeta,\cdot\right)  :=\operatorname*{e}%
\nolimits^{\operatorname*{i}\left(  \operatorname{Im}\zeta\right)
\left\langle \cdot,c\right\rangle }L_{%
%TCIMACRO{\TeXButton{boldmu}{\mbox{\boldmath$ \mu$}}}%
%BeginExpansion
\mbox{\boldmath$ \mu$}%
%EndExpansion
}^{t}\quad\forall t\in\mathcal{T}_{\mathcal{I}},\quad\forall%
%TCIMACRO{\TeXButton{boldmu}{\mbox{\boldmath$ \mu$}}}%
%BeginExpansion
\mbox{\boldmath$ \mu$}%
%EndExpansion
\in\mathbb{N}_{t},\quad\forall c\in\mathcal{D}_{\operatorname*{level}\left(
t\right)  }%
\]
and the expansion coefficients for $b=\left(  t,s\right)  $ by%
\begin{equation}
\gamma_{%
%TCIMACRO{\TeXButton{boldmu}{\mbox{\boldmath$ \mu$}}}%
%BeginExpansion
\mbox{\boldmath$ \mu$}%
%EndExpansion
,%
%TCIMACRO{\TeXButton{boldnu}{\mbox{\boldmath$ \nu$}}}%
%BeginExpansion
\mbox{\boldmath$ \nu$}%
%EndExpansion
,c\left(  b\right)  }^{\left(  t,s\right)  }\left(  \zeta\right)
:=G_{c\left(  b\right)  }\left(  \zeta,\xi_{%
%TCIMACRO{\TeXButton{boldmu}{\mbox{\boldmath$ \mu$}}}%
%BeginExpansion
\mbox{\boldmath$ \mu$}%
%EndExpansion
,t}-\xi_{%
%TCIMACRO{\TeXButton{boldnu}{\mbox{\boldmath$ \nu$}}}%
%BeginExpansion
\mbox{\boldmath$ \nu$}%
%EndExpansion
,s}\right)  \quad\forall b=\left(  t,s\right)  \in\mathcal{P}%
_{\operatorname*{far}},\text{\quad}\forall%
%TCIMACRO{\TeXButton{boldmu}{\mbox{\boldmath$ \mu$}}}%
%BeginExpansion
\mbox{\boldmath$ \mu$}%
%EndExpansion
,%
%TCIMACRO{\TeXButton{boldnu}{\mbox{\boldmath$ \nu$}}}%
%BeginExpansion
\mbox{\boldmath$ \nu$}%
%EndExpansion
\in\mathbb{N}_{t}\text{.}\label{expcoeff}%
\end{equation}
Although this approximation will be slightly modified we introduce the (first)
approximate matrix representation of the sesquilinear form $a\left(
\cdot,\cdot\right)  :S\times S\rightarrow\mathbb{C}$ (cf.
(\ref{weakform_discrete}). Let $\varphi,\psi\in S$ denote some boundary
element functions with basis representation%
\begin{equation}
\varphi=\sum_{i=1}^{n}\varphi_{i}b_{i}\quad\text{and\quad}\psi=\sum_{i=1}%
^{n}\psi_{i}b_{i}.\label{defpsiphi}%
\end{equation}
The coefficients are collected in $%
%TCIMACRO{\TeXButton{boldphi}{\mbox{\boldmath$ \phi$}}}%
%BeginExpansion
\mbox{\boldmath$ \phi$}%
%EndExpansion
=\left(  \varphi_{i}\right)  _{i=1}^{n}$ and $%
%TCIMACRO{\TeXButton{boldpsi}{\mbox{\boldmath$ \psi$}}}%
%BeginExpansion
\mbox{\boldmath$ \psi$}%
%EndExpansion
=\left(  \psi_{i}\right)  _{i=1}^{n}$. We define the (sparse) matrix
$\mathbf{K}^{\operatorname*{near}}=\left(  K_{i,j}^{\operatorname*{near}%
}\right)  _{i,j=1}^{n}\in\mathbb{C}^{n\times n}$ by%
\begin{equation}
K_{i,j}^{\operatorname*{near}}:=\left\{
\begin{array}
[c]{ll}%
K_{i,j} & \text{if }\left(  i,j\right)  \in\left(  \hat{t},\hat{s}\right)
\text{ for some block }\left(  t,s\right)  \in\mathcal{P}%
_{\operatorname*{near}}\text{,}\\
0 & \text{otherwise.}%
\end{array}
\right.  .\label{nearfield_matrix}%
\end{equation}
Then%
\[
a_{\zeta}\left(  \varphi,\psi\right)  \approx\left\langle \mathbf{K}%
^{\operatorname*{near}}%
%TCIMACRO{\TeXButton{boldphi}{\mbox{\boldmath$ \phi$}}}%
%BeginExpansion
\mbox{\boldmath$ \phi$}%
%EndExpansion
,\overline{%
%TCIMACRO{\TeXButton{boldpsi}{\mbox{\boldmath$ \psi$}}}%
%BeginExpansion
\mbox{\boldmath$ \psi$}%
%EndExpansion
}\right\rangle +\sum_{b=\left(  t,s\right)  \in\mathcal{P}%
_{\operatorname*{far}}}\sum_{%
%TCIMACRO{\TeXButton{boldmu}{\mbox{\boldmath$ \mu$}}}%
%BeginExpansion
\mbox{\boldmath$ \mu$}%
%EndExpansion
,%
%TCIMACRO{\TeXButton{boldnu}{\mbox{\boldmath$ \nu$}}}%
%BeginExpansion
\mbox{\boldmath$ \nu$}%
%EndExpansion
\in\mathbb{N}_{t}}\gamma_{%
%TCIMACRO{\TeXButton{boldmu}{\mbox{\boldmath$ \mu$}}}%
%BeginExpansion
\mbox{\boldmath$ \mu$}%
%EndExpansion
,%
%TCIMACRO{\TeXButton{boldnu}{\mbox{\boldmath$ \nu$}}}%
%BeginExpansion
\mbox{\boldmath$ \nu$}%
%EndExpansion
,c}^{b}\left(  \zeta\right)  \tilde{J}_{%
%TCIMACRO{\TeXButton{boldmu}{\mbox{\boldmath$ \mu$}}}%
%BeginExpansion
\mbox{\boldmath$ \mu$}%
%EndExpansion
,c}^{t}\left(  \zeta,\psi\right)  \overline{\tilde{J}_{%
%TCIMACRO{\TeXButton{boldnu}{\mbox{\boldmath$ \nu$}}}%
%BeginExpansion
\mbox{\boldmath$ \nu$}%
%EndExpansion
,c}^{s}\left(  \zeta,\varphi\right)  }%
\]
where $c=c\left(  b\right)  $ and the \textit{farfield coefficients are given
by}%
\[
\tilde{J}_{%
%TCIMACRO{\TeXButton{boldmu}{\mbox{\boldmath$ \mu$}}}%
%BeginExpansion
\mbox{\boldmath$ \mu$}%
%EndExpansion
,c}^{t}\left(  \zeta,\psi\right)  :=\sum_{i\in\hat{t}}\overline{\psi_{i}}%
\int_{\Gamma}\tilde{\Phi}_{%
%TCIMACRO{\TeXButton{boldmu}{\mbox{\boldmath$ \mu$}}}%
%BeginExpansion
\mbox{\boldmath$ \mu$}%
%EndExpansion
,c}^{t}\left(  \zeta,x\right)  b_{i}\left(  x\right)  d\Gamma_{x}.
\]

Since the expansion orders $m_{\ell}$ are monotonously increasing we can
express a Lagrange basis $L_{%
%TCIMACRO{\TeXButton{boldmu}{\mbox{\boldmath$ \mu$}}}%
%BeginExpansion
\mbox{\boldmath$ \mu$}%
%EndExpansion
}^{t}$ via the Lagrange basis on $t^{\prime}\in\operatorname*{sons}\left(
t\right)  :$%
\begin{equation}
L_{%
%TCIMACRO{\TeXButton{boldmu}{\mbox{\boldmath$ \mu$}}}%
%BeginExpansion
\mbox{\boldmath$ \mu$}%
%EndExpansion
}^{t}=\sum_{%
%TCIMACRO{\TeXButton{boldnu}{\mbox{\boldmath$ \nu$}}}%
%BeginExpansion
\mbox{\boldmath$ \nu$}%
%EndExpansion
\in\mathbb{N}_{t^{\prime}}}q_{t^{\prime},%
%TCIMACRO{\TeXButton{boldmu}{\mbox{\boldmath$ \mu$}}}%
%BeginExpansion
\mbox{\boldmath$ \mu$}%
%EndExpansion
,%
%TCIMACRO{\TeXButton{boldnu}{\mbox{\boldmath$ \nu$}}}%
%BeginExpansion
\mbox{\boldmath$ \nu$}%
%EndExpansion
}L_{%
%TCIMACRO{\TeXButton{boldnu}{\mbox{\boldmath$ \nu$}}}%
%BeginExpansion
\mbox{\boldmath$ \nu$}%
%EndExpansion
}^{t^{\prime}} \label{defq}%
\end{equation}
with the \textit{transfer coefficients} $q_{t^{\prime},%
%TCIMACRO{\TeXButton{boldmu}{\mbox{\boldmath$ \mu$}}}%
%BeginExpansion
\mbox{\boldmath$ \mu$}%
%EndExpansion
,%
%TCIMACRO{\TeXButton{boldnu}{\mbox{\boldmath$ \nu$}}}%
%BeginExpansion
\mbox{\boldmath$ \nu$}%
%EndExpansion
}=L_{%
%TCIMACRO{\TeXButton{boldmu}{\mbox{\boldmath$ \mu$}}}%
%BeginExpansion
\mbox{\boldmath$ \mu$}%
%EndExpansion
}^{t}\left(  \xi_{%
%TCIMACRO{\TeXButton{boldnu}{\mbox{\boldmath$ \nu$}}}%
%BeginExpansion
\mbox{\boldmath$ \nu$}%
%EndExpansion
,t^{\prime}}\right)  $. We cannot expect that a direction $c\in\mathcal{D}%
_{\ell}$ is also contained in the set $\mathcal{D}_{\ell+1}$, hence we
\textit{assign} to $c$ the direction $c^{\prime}=\operatorname*{sd}\left(
c\right)  \in\mathcal{D}_{\ell+1}$ which has a minimal Euclidean distance.
This leads to the recursive definition of the final expansion functions%
\[
\Phi_{%
%TCIMACRO{\TeXButton{boldmu}{\mbox{\boldmath$ \mu$}}}%
%BeginExpansion
\mbox{\boldmath$ \mu$}%
%EndExpansion
,c}^{t}\left(  \zeta,\cdot\right)  :=\operatorname*{e}%
\nolimits^{-\operatorname*{i}\left(  \operatorname{Im}\zeta\right)
\left\langle \cdot,c\right\rangle }L_{%
%TCIMACRO{\TeXButton{boldmu}{\mbox{\boldmath$ \mu$}}}%
%BeginExpansion
\mbox{\boldmath$ \mu$}%
%EndExpansion
}^{t}\quad\forall t\in\mathcal{L}_{\mathcal{I}},\quad\forall%
%TCIMACRO{\TeXButton{boldmu}{\mbox{\boldmath$ \mu$}}}%
%BeginExpansion
\mbox{\boldmath$ \mu$}%
%EndExpansion
\in\mathbb{N}_{t}\text{,\quad}\forall c\in\mathcal{D}_{\operatorname*{level}%
\left(  t\right)  }%
\]
and for $t\in\mathcal{T}_{\mathcal{I}}\backslash\mathcal{L}_{\mathcal{I}}$ and
$t^{\prime}\in\operatorname*{sons}\left(  t\right)  $ we set $c^{\prime
}:=\operatorname*{sd}\left(  c\right)  $ and%
\begin{equation}
\Phi_{%
%TCIMACRO{\TeXButton{boldmu}{\mbox{\boldmath$ \mu$}}}%
%BeginExpansion
\mbox{\boldmath$ \mu$}%
%EndExpansion
,c}^{t}\left(  \zeta,\cdot\right)  :=\operatorname*{e}%
\nolimits^{\operatorname*{i}\left(  \operatorname{Im}\zeta\right)
\left\langle \cdot,c^{\prime}\right\rangle }\sum_{%
%TCIMACRO{\TeXButton{boldnu}{\mbox{\boldmath$ \nu$}}}%
%BeginExpansion
\mbox{\boldmath$ \nu$}%
%EndExpansion
\in\mathbb{N}_{t^{\prime}}}q_{t^{\prime},%
%TCIMACRO{\TeXButton{boldmu}{\mbox{\boldmath$ \mu$}}}%
%BeginExpansion
\mbox{\boldmath$ \mu$}%
%EndExpansion
,%
%TCIMACRO{\TeXButton{boldnu}{\mbox{\boldmath$ \nu$}}}%
%BeginExpansion
\mbox{\boldmath$ \nu$}%
%EndExpansion
}L_{%
%TCIMACRO{\TeXButton{boldnu}{\mbox{\boldmath$ \nu$}}}%
%BeginExpansion
\mbox{\boldmath$ \nu$}%
%EndExpansion
}^{t^{\prime}}. \label{rekphi}%
\end{equation}
This, in turn, motivates the definition of the farfield coefficients
corresponding to $\Phi_{%
%TCIMACRO{\TeXButton{boldmu}{\mbox{\boldmath$ \mu$}}}%
%BeginExpansion
\mbox{\boldmath$ \mu$}%
%EndExpansion
,c}^{t}$ by%
\[
J_{%
%TCIMACRO{\TeXButton{boldmu}{\mbox{\boldmath$ \mu$}}}%
%BeginExpansion
\mbox{\boldmath$ \mu$}%
%EndExpansion
,c}^{t}\left(  \zeta,\psi\right)  :=\sum_{i\in\hat{t}}\overline{\psi_{i}}%
\int_{\Gamma}\Phi_{%
%TCIMACRO{\TeXButton{boldmu}{\mbox{\boldmath$ \mu$}}}%
%BeginExpansion
\mbox{\boldmath$ \mu$}%
%EndExpansion
,c}^{t}\left(  \zeta,x\right)  b_{i}\left(  x\right)  d\Gamma_{x}.
\]
The relation (\ref{rekphi}) allows for an hierarchical computation of these
coefficients. First we compute the basis farfield coefficients%
\begin{equation}
J_{%
%TCIMACRO{\TeXButton{boldmu}{\mbox{\boldmath$ \mu$}}}%
%BeginExpansion
\mbox{\boldmath$ \mu$}%
%EndExpansion
,c}^{t}\left(  \zeta,b_{j}\right)  :=\int_{\Gamma}\Phi_{%
%TCIMACRO{\TeXButton{boldmu}{\mbox{\boldmath$ \mu$}}}%
%BeginExpansion
\mbox{\boldmath$ \mu$}%
%EndExpansion
,c}^{t}\left(  \zeta,x\right)  b_{j}\left(  x\right)  d\Gamma_{x}\qquad\forall
t\in\mathcal{L}_{\mathcal{I}}\quad\forall j\in\hat{t}\quad\forall
c\in\mathcal{D}_{\operatorname*{level}\left(  t\right)  }\quad\forall%
%TCIMACRO{\TeXButton{boldmu}{\mbox{\boldmath$ \mu$}}}%
%BeginExpansion
\mbox{\boldmath$ \mu$}%
%EndExpansion
\in\mathbb{N}_{t}. \label{basisffc}%
\end{equation}
Then, for a boundary element function $\psi$ as in (\ref{defpsiphi}) we
determine%
\[
J_{%
%TCIMACRO{\TeXButton{boldmu}{\mbox{\boldmath$ \mu$}}}%
%BeginExpansion
\mbox{\boldmath$ \mu$}%
%EndExpansion
,c}^{t}\left(  \zeta,\psi\right)  =\sum_{i\in\hat{t}}\overline{\psi_{i}}J_{%
%TCIMACRO{\TeXButton{boldmu}{\mbox{\boldmath$ \mu$}}}%
%BeginExpansion
\mbox{\boldmath$ \mu$}%
%EndExpansion
,c}^{t}\left(  \zeta,b_{i}\right)  \qquad\forall t\in\mathcal{L}_{\mathcal{I}%
}\quad\forall c\in\mathcal{D}_{\operatorname*{level}\left(  t\right)  }%
\quad\forall%
%TCIMACRO{\TeXButton{boldmu}{\mbox{\boldmath$ \mu$}}}%
%BeginExpansion
\mbox{\boldmath$ \mu$}%
%EndExpansion
\in\mathbb{N}_{t}%
\]
and recursively%
\[
J_{%
%TCIMACRO{\TeXButton{boldmu}{\mbox{\boldmath$ \mu$}}}%
%BeginExpansion
\mbox{\boldmath$ \mu$}%
%EndExpansion
,c}^{t}\left(  \zeta,\psi\right)  :=\sum_{t^{\prime}\in\operatorname*{sons}%
\left(  t\right)  }\sum_{%
%TCIMACRO{\TeXButton{boldnu}{\mbox{\boldmath$ \nu$}}}%
%BeginExpansion
\mbox{\boldmath$ \nu$}%
%EndExpansion
\in\mathbb{N}_{t^{\prime}}}q_{t^{\prime},%
%TCIMACRO{\TeXButton{boldmu}{\mbox{\boldmath$ \mu$}}}%
%BeginExpansion
\mbox{\boldmath$ \mu$}%
%EndExpansion
,%
%TCIMACRO{\TeXButton{boldnu}{\mbox{\boldmath$ \nu$}}}%
%BeginExpansion
\mbox{\boldmath$ \nu$}%
%EndExpansion
}J_{%
%TCIMACRO{\TeXButton{boldmu}{\mbox{\boldmath$ \mu$}}}%
%BeginExpansion
\mbox{\boldmath$ \mu$}%
%EndExpansion
,c^{\prime}}^{t^{\prime}}\left(  \zeta,\psi\right)  \qquad\forall
t\in\mathcal{T}_{\mathcal{I}}\backslash\mathcal{L}_{\mathcal{I}}\quad\forall
c\in\mathcal{D}_{\operatorname*{level}\left(  t\right)  }\quad\forall%
%TCIMACRO{\TeXButton{boldmu}{\mbox{\boldmath$ \mu$}}}%
%BeginExpansion
\mbox{\boldmath$ \mu$}%
%EndExpansion
\in\mathbb{N}_{t}%
\]
by using the hierarchical tree structure.

Once, the farfield coefficients are $J_{%
%TCIMACRO{\TeXButton{boldmu}{\mbox{\boldmath$ \mu$}}}%
%BeginExpansion
\mbox{\boldmath$ \mu$}%
%EndExpansion
,c}^{t}$ computed, the final approximation of the sesquilinear form $a_{\zeta
}\left(  \cdot,\cdot\right)  $ can be evaluated%
\begin{align}
a_{\zeta}\left(  \varphi,\psi\right)   &  \approx a_{\zeta}^{\mathcal{DH}^{2}%
}\left(  \varphi,\psi\right)  :=\left\langle \mathbf{K}^{\operatorname*{near}}%
%TCIMACRO{\TeXButton{boldphi}{\mbox{\boldmath$ \phi$}}}%
%BeginExpansion
\mbox{\boldmath$ \phi$}%
%EndExpansion
,%
%TCIMACRO{\TeXButton{boldpsi}{\mbox{\boldmath$ \psi$}}}%
%BeginExpansion
\mbox{\boldmath$ \psi$}%
%EndExpansion
\right\rangle \label{h2approxazeta}\\
&  +\sum_{b=\left(  t,s\right)  \in\mathcal{P}_{\operatorname*{far}}}\sum_{%
%TCIMACRO{\TeXButton{boldmu}{\mbox{\boldmath$ \mu$}}}%
%BeginExpansion
\mbox{\boldmath$ \mu$}%
%EndExpansion
,%
%TCIMACRO{\TeXButton{boldnu}{\mbox{\boldmath$ \nu$}}}%
%BeginExpansion
\mbox{\boldmath$ \nu$}%
%EndExpansion
\in\mathbb{N}_{t}}\gamma_{%
%TCIMACRO{\TeXButton{boldmu}{\mbox{\boldmath$ \mu$}}}%
%BeginExpansion
\mbox{\boldmath$ \mu$}%
%EndExpansion
,%
%TCIMACRO{\TeXButton{boldnu}{\mbox{\boldmath$ \nu$}}}%
%BeginExpansion
\mbox{\boldmath$ \nu$}%
%EndExpansion
,c\left(  b\right)  }^{b}\left(  \zeta\right)  J_{%
%TCIMACRO{\TeXButton{boldmu}{\mbox{\boldmath$ \mu$}}}%
%BeginExpansion
\mbox{\boldmath$ \mu$}%
%EndExpansion
,c\left(  b\right)  }^{t}\left(  \zeta,\psi\right)  \overline{J_{%
%TCIMACRO{\TeXButton{boldnu}{\mbox{\boldmath$ \nu$}}}%
%BeginExpansion
\mbox{\boldmath$ \nu$}%
%EndExpansion
,c\left(  b\right)  }^{s}\left(  \zeta,\varphi\right)  }.\nonumber
\end{align}
The algorithmic formulation of an approximate matrix vector multiplication,
i.e., the computation of $\left(  a_{\zeta}\left(  \varphi,b_{i}\right)
\right)  _{i=1}^{n}$ can be derived from (\ref{h2approxazeta}) and the details
are in the literature, e.g., in \cite{SauterDiss}, \cite{Sauter00},
\cite{SauterSchwab2010}, and for our concrete application, e.g., in
\cite{Boerm2015}.

\begin{remark}
We will prove in Sections \ref{SecAna} and \ref{SecComplexity} that the
compression algorithm presented in this section results in a sparse
$\mathcal{DH}^{2}$-matrix approximation and analyze how to choose the control
parameters in order to satisfy a prescribed accuracy for this perturbation.
However, numerical experiments show that the rank of this approximation may be
larger than necessary. In \cite{BoermMelenk}, \cite{borm2018hybrid} a
\emph{recompression algorithm} is presented for the pure Helmholtz problem
($\zeta\in\operatorname*{i}\mathbb{R}$) which further compresses an already
sparse $\mathcal{DH}^{2}$-matrix. We do not elaborate this issue here since
the recompression algorithm in \cite{BoermMelenk}, \cite{borm2018hybrid} can
be applied verbatim to the case of general complex frequencies and results in
nearly optimal storage requirements. This recompression algorithm on top of
our $\mathcal{DH}^{2}$-matrix approximation will be employed for our numerical
experiments in Section \ref{SecNumExp}.
\end{remark}

\section{Analysis\label{SecAna}}

In this section, we will investigate the accuracy of the directional
$\mathcal{H}^{2}$ approximation for the acoustic single layer potential for
general complex frequencies $z\in\mathbb{C}_{>0}$ by generalizing the results
in \cite{BoermMelenk}.

The key role is played by derivative-free interpolation estimates which go
back to \cite{Dr}. Let $\left(  t,s\right)  $ denote an admissible block. For
$x\in B_{t}$ and $y\in B_{s}$, let $z=y-x$ and $r=\left\Vert z\right\Vert $.
The kernel function of the acoustic single layer potential for the complex
frequency $\zeta\in\mathbb{C}_{>0}$, is given by%
\[
G\left(  \zeta,z\right)  :=\frac{\operatorname*{e}^{-\zeta r}}{4\pi
r}=\operatorname*{e}\nolimits^{-\operatorname*{i}\left(  \operatorname*{Im}%
\zeta\right)  \left\langle z,c\right\rangle }G_{c}\left(  z\right)
\]
with%
\[
G_{c}\left(  z\right)  :=\operatorname*{e}\nolimits^{-\left(
\operatorname{Re}\zeta\right)  r}\frac{\operatorname*{e}%
\nolimits^{-\operatorname*{i}\left(  \operatorname*{Im}\zeta\right)  \left(
r-\left\langle z,c\right\rangle \right)  }}{4\pi r}%
\]
for some unit vector $c\in\mathbb{R}^{3}$, $\left\Vert c\right\Vert =1$. We
approximate this function on $B_{t}\times B_{s}$ by%
\[
\tilde{G}_{t,s}\left(  z\right)  :=\operatorname*{e}%
\nolimits^{-\operatorname*{i}\left(  \operatorname{Im}\zeta\right)
\left\langle z,c\right\rangle }\mathfrak{I}_{t\times s}\left(  G_{c}\right)
,
\]
where $\mathfrak{I}_{t\times s}$ denote the tensor \v{C}eby\v{s}ev
interpolation on $B_{t}\times B_{s}$ with polynomials of maximal degree $m$.

For the error analysis, we modify the theory as in \cite{BoermMelenk} and
present the relevant statements in the following.

\begin{lemma}
\label{Lempest}Let $t$, $s$, and $c$ satisfy the conditions (\ref{adm_angle}),
(\ref{adm_paranew}). Let $d,p\in\mathbb{R}^{3}$ be vectors satisfying%
%TCIMACRO{\TeXButton{30}{\begin{subequations}
%\label{30}
%\end{subequations}} }%
%BeginExpansion
\begin{subequations}
\label{30}
\end{subequations}
%EndExpansion
\begin{align}
\left\Vert p\right\Vert  &  \leq\frac{\max\left\{  \operatorname*{diam}\left(
B_{t}\right)  ,\operatorname*{diam}\left(  B_{s}\right)  \right\}  }{2},\tag{%
%TCIMACRO{\TeXButton{30}{\ref{30}}}%
%BeginExpansion
\ref{30}%
%EndExpansion
a}\label{30a}\\
d-\tau p  &  \in B_{t}-B_{s}=\left\{  x-y:\left(  x,y\right)  \in B_{t}\times
B_{s}\right\}  \quad\forall\tau\in\left[  -1,1\right]  . \tag{%
%TCIMACRO{\TeXButton{30}{\ref{30}}}%
%BeginExpansion
\ref{30}%
%EndExpansion
b}\label{30b}%
\end{align}
Then, we have%
\[
\left\Vert \frac{d-\tau p}{\left\Vert d-\tau p\right\Vert }-c\right\Vert
\leq\frac{\eta_{1}+\max\left\{  \eta_{2},\eta_{3}\left(  \operatorname{Re}%
\zeta\right)  \operatorname*{dist}\left(  B_{t},B_{s}\right)  \right\}
}{\left\vert \operatorname{Im}\zeta\right\vert \max\left\{
\operatorname*{diam}\left(  B_{t}\right)  ,\operatorname*{diam}\left(
B_{s}\right)  \right\}  }\quad\forall\tau\in\left[  -1,1\right]  .
\]

\end{lemma}

%

%TCIMACRO{\TeXButton{Proof}{\proof}}%
%BeginExpansion
\proof
%EndExpansion
We only sketch the minor modifications in the proof of \cite[Lemma
3.9]{BoermMelenk} for our modified admissibility condition (cf.
(\ref{adm_paranew})). It holds%
\begin{align*}
\left\Vert d-\tau p\right\Vert  &  \geq\operatorname*{dist}\left(  B_{t}%
,B_{s}\right)  \geq\frac{\left\vert \operatorname{Im}\zeta\right\vert q^{2}%
}{\check{\eta}},\\
\left\Vert M_{t}-M_{s}\right\Vert  &  \geq\operatorname*{dist}\left(
B_{t},B_{s}\right)  \geq\frac{\left\vert \operatorname{Im}\zeta\right\vert
q^{2}}{\check{\eta}}.
\end{align*}
with $q:=\max\left\{  \operatorname*{diam}\left(  B_{t}\right)
,\operatorname*{diam}\left(  B_{s}\right)  \right\}  $ and $\check{\eta}%
:=\max\left\{  \eta_{2},\eta_{3}\left(  \operatorname{Re}\zeta\right)
\operatorname*{dist}\left(  B_{t},B_{s}\right)  \right\}  $. Hence,%
\[
\left\Vert \frac{d-\tau p}{\left\Vert d-\tau p\right\Vert }-\frac{M_{t}-M_{s}%
}{\left\Vert M_{t}-M_{s}\right\Vert }\right\Vert \leq\frac{\check{\eta}%
}{\left\vert \operatorname{Im}\zeta\right\vert q}%
\]
and%
\begin{align*}
\left\Vert \frac{d-\tau p}{\left\Vert d-\tau p\right\Vert }-c\right\Vert  &
\leq\left\Vert \frac{d-\tau p}{\left\Vert d-\tau p\right\Vert }-\frac
{M_{t}-M_{s}}{\left\Vert M_{t}-M_{s}\right\Vert }\right\Vert +\left\Vert
\frac{M_{t}-M_{s}}{\left\Vert M_{t}-M_{s}\right\Vert }-c\right\Vert \\
&  \leq\frac{\check{\eta}}{\left\vert \operatorname{Im}\zeta\right\vert
q}+\frac{\eta_{1}}{\left\vert \operatorname{Im}\zeta\right\vert q}.
\end{align*}%
%TCIMACRO{\TeXButton{End Proof}{\endproof}}%
%BeginExpansion
\endproof
%EndExpansion

To formulate the main theorem for the interpolation error, we introduce first
some constants. Let $\eta_{1}$, $\eta_{2}$, $\eta_{3}$ denote the positive
control parameters for the admissibility conditions (\ref{adm}). We assume
that $0<\eta_{3}<1$ holds and set $\sigma:=1-\eta_{3}>0$. The Lebesgue
constant for the univariate \v{C}eby\v{s}ev interpolation is denoted by
$\Lambda_{m}$ and we recall the well-known estimate $\Lambda_{m}\leq\frac
{2}{\pi}\log\left(  m+1\right)  +1$. Let%
\begin{equation}
\rho_{0}:=1+\widehat{\beta}\quad\text{and\quad}\widehat{\beta}:=\min\left\{
1,\left(  \sqrt{\frac{3}{2}}-1\right)  \frac{2}{\eta_{2}},\frac{2\left(
1-\eta_{3}\right)  }{\eta_{2}^{2}\left(  2\sqrt{6}+5\right)  }\right\}  .
\label{defbetahat}%
\end{equation}
We set%
\begin{align}
\alpha &  :=\left(  \sqrt{\widehat{\beta}^{2}+1}+\widehat{\beta}\right)
/\left(  \widehat{\beta}+1\right) \label{defalpha}\\
C_{1}  &  :=\operatorname*{e}\nolimits^{\eta_{1}}\sup\left\{  \frac{8\left(
\Lambda_{m}+1\right)  }{\left(  \rho_{0}-1\right)  \alpha^{m/2}}%
:m\in\mathbb{N}\right\}  ,\label{defC1}\\
C_{0}  &  :=\sup\left\{  6\frac{\Lambda_{m}^{5}}{\alpha^{m/2}}C_{1}%
:m\in\mathbb{N}\right\}
\end{align}
and observe that $C_{0}$, $C_{1}$ are bounded since the Lebesgue constant
grows only logarithmically in $m$ and it is easy to see that $\alpha>1$.

\begin{theorem}
\label{TheoErrLoc}Let $c\in\mathbb{R}^{3}$ and let the block $b=\left(
t,s\right)  $ satisfy the $%
%TCIMACRO{\TeXButton{boldeta}{\mbox{\boldmath$ \eta$}}}%
%BeginExpansion
\mbox{\boldmath$ \eta$}%
%EndExpansion
$-admissibility conditions (\ref{adm}) for some $\eta_{3}\in\left(
0,1\right)  $. Then,%
\[
\left\Vert G-\tilde{G}_{t,s}\right\Vert _{\infty,t\times s}\leq\frac
{C_{0}\operatorname*{e}^{-\sigma\left(  \operatorname{Re}\zeta\right)
\operatorname*{dist}\left(  B_{t},B_{s}\right)  }}{4\pi\operatorname*{dist}%
\left(  B_{t},B_{s}\right)  }\rho_{0}^{-m}\quad\text{for }\sigma:=\frac
{1-\eta_{3}}{2}%
\]
where $\rho_{0}>1$ depends on $\eta_{2}$ and $\eta_{3}$.
\end{theorem}

%

%TCIMACRO{\TeXButton{Proof}{\proof}}%
%BeginExpansion
\proof
%EndExpansion
We introduce the function $G_{\operatorname*{dp}}:\left[  -1,1\right]
\rightarrow\mathbb{C}$ by%
\[
G_{\operatorname*{dp}}\left(  x\right)  :=\operatorname*{e}\nolimits^{-\left(
\operatorname{Re}\zeta\right)  \left\Vert d-xp\right\Vert }\frac{\exp\left(
-\operatorname*{i}\left(  \operatorname{Im}\zeta\right)  \left(  \left\Vert
d-xp\right\Vert -\left\langle d-xp,c\right\rangle \right)  \right)  }%
{4\pi\left\Vert d-xp\right\Vert }.
\]
and first prove an error estimate for the univariate \v{C}eby\v{s}ev
interpolation of this function. Lemma \ref{Lempest} leads to the estimate%
\[
\xi:=\max\left\{  \left\Vert \frac{d-\gamma p}{\left\Vert d-\gamma
p\right\Vert }-c\right\Vert :\gamma\in\left[  -1,1\right]  \right\}  \leq
\frac{\eta_{1}+\max\left\{  \eta_{2},\eta_{3}\left(  \operatorname{Re}%
\zeta\right)  \operatorname*{dist}\left(  B_{t},B_{s}\right)  \right\}
}{\left\vert \operatorname{Im}\zeta\right\vert \max\left\{
\operatorname*{diam}\left(  B_{t}\right)  ,\operatorname*{diam}\left(
B_{s}\right)  \right\}  }.
\]
Note that the admissibility conditions imply that (\ref{30}) holds. Hence,
(\ref{30b}) yields%
\[
\delta:=\inf\left\{  \left\Vert d-\gamma p\right\Vert :\gamma\in\left[
-1,1\right]  \right\}  \geq\operatorname*{dist}\left(  B_{t},B_{s}\right)  .
\]
The combination with (\ref{30a}) leads to%
\begin{equation}
\lambda:=\delta/\left\Vert p\right\Vert \geq\frac{2\operatorname*{dist}\left(
B_{t},B_{s}\right)  }{\max\left\{  \operatorname*{diam}\left(  B_{t}\right)
,\operatorname*{diam}B_{s}\right\}  }\geq2/\eta_{2}>0\quad\text{i.e.,\quad
}1/\lambda\leq\eta_{2}/2. \label{lambdaestp}%
\end{equation}
To estimate the interpolation error for the function $G_{\operatorname*{dp}}$
we first derive a bound for the modulus of $G_{\operatorname*{dp}}$ in a
complex neighborhood of $\left[  -1,1\right]  $. We set%
\[
\beta:=\min\left\{  1,\left(  \sqrt{\frac{3}{2}}-1\right)  \lambda
,\frac{\lambda\left(  1-\eta_{3}\right)  }{\eta_{2}\left(  2\sqrt{6}+5\right)
}\right\}  \overset{\lambda\geq2/\eta_{2}}{\geq}\widehat{\beta}%
\]
with $\widehat{\beta}$ as in (\ref{defbetahat}) and define $U_{\beta
}:=\left\{  z\in\mathbb{C}\mid\operatorname*{dist}\left(  z,\left[
-1,1\right]  \right)  \leq\beta\right\}  $. The unique analytic continuation
of the square root function $\sqrt{\cdot}:\mathbb{R}_{>0}\rightarrow
\mathbb{R}_{>0}$ to $\mathbb{C}\backslash\mathbb{R}_{\leq0}$ is given by%
\[
\sqrt{z}=\sqrt{\left\vert z\right\vert }\frac{z+\left\vert z\right\vert
}{\left\vert z+\left\vert z\right\vert \right\vert }\qquad\forall
z\in\mathbb{C}\backslash\mathbb{R}_{\leq0}.
\]
The analytic continuation of the function $x\rightarrow\left\Vert
d-xp\right\Vert ^{2}$ is then denoted by%
\[
\mathfrak{n}_{\operatorname*{dp}}\left(  z\right)  :=\sqrt{\left\langle
d-zp,d-zp\right\rangle }.
\]
The modulus of $G_{\operatorname*{dp}}$ can be estimated by%
\[
\sup_{z\in U_{\beta}}\left\vert G_{\operatorname*{dp}}\left(  z\right)
\right\vert \leq w\left(  \beta\right)  \chi_{+}\left(  \beta\right)  \chi
_{-}\left(  \beta\right)  .
\]
with $w\left(  \beta\right)  :=\sup_{z\in U_{\beta}}\operatorname*{e}%
\nolimits^{-\left(  \operatorname{Re}\zeta\right)  \operatorname{Re}\left(
\mathfrak{n}_{\operatorname*{dp}}\left(  z\right)  \right)  }$ and%
\[
\chi_{+}\left(  \beta\right)  :=\sup_{z\in U_{\beta}}\exp\left(
-\operatorname*{i}\left(  \operatorname{Im}\zeta\right)  \left(
\mathfrak{n}_{\operatorname*{dp}}\left(  z\right)  -\left\langle
d-zp,c\right\rangle \right)  \right)  ,\quad\chi_{-}\left(  \beta\right)
:=\sup_{z\in U_{\beta}}\frac{1}{4\pi\mathfrak{n}_{\operatorname*{dp}}\left(
z\right)  }.
\]
For $\chi_{+}\left(  \beta\right)  $ we obtain $\chi_{+}\left(  \beta\right)
\leq\exp\left(  \gamma\left(  \beta\right)  \right)  $, where the exponent
$\gamma\left(  \beta\right)  $ can be estimated by (cf. the proof of Lemma 3.8
in \cite{BoermMelenk})%
\begin{align*}
\gamma\left(  \beta\right)   &  =\left\vert \operatorname{Im}\zeta\right\vert
\left\Vert p\right\Vert \beta\left(  \xi+\frac{\beta}{2\left(  \lambda
-\beta\right)  }\right) \\
&  \leq\left\vert \operatorname{Im}\zeta\right\vert \left\Vert p\right\Vert
\beta\frac{\eta_{1}+\max\left\{  \eta_{2},\eta_{3}\left(  \operatorname{Re}%
\zeta\right)  \operatorname*{dist}\left(  B_{t},B_{s}\right)  \right\}
}{\left\vert \operatorname{Im}\zeta\right\vert \max\left\{
\operatorname*{diam}\left(  B_{t}\right)  ,\operatorname*{diam}\left(
B_{s}\right)  \right\}  }+\frac{\left\vert \operatorname{Im}\zeta\right\vert
\left\Vert p\right\Vert \beta^{2}}{2\left(  \lambda-\beta\right)  }\\
&  \overset{\text{(\ref{30a})}}{\leq}\frac{\beta}{2}\left(  \eta_{1}%
+\max\left\{  \eta_{2},\eta_{3}\left(  \operatorname{Re}\zeta\right)
\operatorname*{dist}\left(  B_{t},B_{s}\right)  \right\}  \right)
+\frac{\left\vert \operatorname{Im}\zeta\right\vert \left\Vert p\right\Vert
\beta^{2}}{2\lambda\left(  1-\beta/\lambda\right)  }.
\end{align*}
We obtain with our new parabolic admissibility condition (\ref{adm_paranew})%
\begin{align*}
\frac{\left\vert \operatorname{Im}\zeta\right\vert \left\Vert p\right\Vert
\beta^{2}}{2\lambda\left(  1-\beta/\lambda\right)  }  &  \overset
{\text{(\ref{lambdaestp})}}{\leq}\frac{\left\vert \operatorname{Im}%
\zeta\right\vert \left\Vert p\right\Vert \beta^{2}}{4\left(  1-\beta
/\lambda\right)  }\frac{\max\left\{  \operatorname*{diam}\left(  B_{t}\right)
,\operatorname*{diam}B_{s}\right\}  }{\operatorname*{dist}\left(  B_{t}%
,B_{s}\right)  }\\
&  \overset{\text{(\ref{30a})}}{\leq}\frac{\left\vert \operatorname{Im}%
\zeta\right\vert \max\left\{  \operatorname*{diam}^{2}\left(  B_{t}\right)
,\operatorname*{diam}^{2}B_{s}\right\}  \beta^{2}}{8\left(  1-\beta
/\lambda\right)  \operatorname*{dist}\left(  B_{t},B_{s}\right)  }\\
&  \overset{\text{(\ref{adm_paranew})}}{\leq}\frac{\beta}{2}\max\left\{
\eta_{2},\eta_{3}\left(  \operatorname{Re}\zeta\right)  \operatorname*{dist}%
\left(  B_{t},B_{s}\right)  \right\}  \frac{\beta}{4\left(  1-\beta
/\lambda\right)  }.
\end{align*}
Since $\beta\leq\min\left\{  1,\frac{3}{4}\lambda\right\}  $, we have derived%
\begin{align*}
\gamma\left(  \beta\right)   &  \leq\frac{\beta}{2}\left(  \eta_{1}%
+\max\left\{  \eta_{2},\eta_{3}\left(  \operatorname{Re}\zeta\right)
\operatorname*{dist}\left(  B_{t},B_{s}\right)  \right\}  \right)
+\frac{\beta}{2}\left(  \max\left\{  \eta_{2},\eta_{3}\left(
\operatorname{Re}\zeta\right)  \operatorname*{dist}\left(  B_{t},B_{s}\right)
\right\}  \right) \\
&  \leq\eta_{1}+\max\left\{  \eta_{2},\eta_{3}\left(  \operatorname{Re}%
\zeta\right)  \operatorname*{dist}\left(  B_{t},B_{s}\right)  \right\}  .
\end{align*}
The estimate%
\[
\chi_{-}\left(  \beta\right)  \leq\frac{1}{4\pi\inf_{z\in U_{\beta}}\left\vert
\mathfrak{n}_{\operatorname*{dp}}\left(  z\right)  \right\vert }\leq\frac
{1}{4\pi\operatorname*{dist}\left(  B_{t},B_{s}\right)  \left(  1-\beta
/\lambda\right)  }%
\]
follows as in the proof of \cite[Lemma 3.6]{BoermMelenk}. We employ
$\beta<\frac{3}{4}\lambda$ so that%
\[
\chi_{-}\left(  \beta\right)  \leq\frac{1}{\pi\operatorname*{dist}\left(
B_{t},B_{s}\right)  }.
\]

Next, we estimate the term $w\left(  \beta\right)  $. For $z\in U_{\beta}$, we
choose $x_{z}\in\left[  -1,1\right]  $ such that%
\[
\min_{x\in\left[  -1,1\right]  }\left\vert z-x\right\vert =\left\vert
z-x_{z}\right\vert .
\]
Hence%
\[
\operatorname{Re}\mathfrak{n}_{\operatorname*{dp}}\left(  z\right)
\geq\left\Vert d-x_{z}p\right\Vert -\left\vert n_{\operatorname*{dp}}\left(
z\right)  -n_{\operatorname*{dp}}\left(  x_{z}\right)  \right\vert .
\]
We set $\psi:=d-zp$ and $\varphi=d-x_{z}p$ so that%
\begin{align*}
n_{\operatorname*{dp}}\left(  z\right)  -n_{\operatorname*{dp}}\left(
x_{z}\right)   &  =\sqrt{\left\langle \psi,\psi\right\rangle }-\left\Vert
\varphi\right\Vert =\left\Vert \psi\right\Vert \frac{\left\langle \psi
,\psi\right\rangle +\left\Vert \psi\right\Vert ^{2}}{\left\vert \left\langle
\psi,\psi\right\rangle +\left\Vert \psi\right\Vert ^{2}\right\vert
}-\left\Vert \varphi\right\Vert \\
&  =\left(  \left\Vert \psi\right\Vert -\left\Vert \varphi\right\Vert \right)
\frac{\left\langle \psi,\psi\right\rangle +\left\Vert \psi\right\Vert ^{2}%
}{\left\vert \left\langle \psi,\psi\right\rangle +\left\Vert \psi\right\Vert
^{2}\right\vert }+\left\Vert \varphi\right\Vert \left(  \frac{\left\langle
\psi,\psi\right\rangle +\left\Vert \psi\right\Vert ^{2}}{\left\vert
\left\langle \psi,\psi\right\rangle +\left\Vert \psi\right\Vert ^{2}%
\right\vert }-1\right)  .
\end{align*}
This leads to the estimate%
\[
\left\vert n_{\operatorname*{dp}}\left(  z\right)  -n_{\operatorname*{dp}%
}\left(  x_{z}\right)  \right\vert \leq\left\vert \left\Vert \psi\right\Vert
-\left\Vert \varphi\right\Vert \right\vert +\left\Vert \varphi\right\Vert
\left\vert \frac{\left\langle \psi,\psi\right\rangle +\left\Vert
\psi\right\Vert ^{2}}{\left\vert \left\langle \psi,\psi\right\rangle
+\left\Vert \psi\right\Vert ^{2}\right\vert }-1\right\vert .
\]
We know $\left\vert z-x_{z}\right\vert \leq\beta$ so that%
\begin{equation}
\left\vert \left\Vert \psi\right\Vert -\left\Vert \varphi\right\Vert
\right\vert \leq\left\Vert p\right\Vert \left\vert z-x_{z}\right\vert
\leq\beta\left\Vert p\right\Vert . \label{normpsimnormvphi}%
\end{equation}
Furthermore, a triangle inequality leads to%
\begin{align*}
&  \left\vert \frac{\left\langle \psi,\psi\right\rangle +\left\Vert
\psi\right\Vert ^{2}}{\left\vert \left\langle \psi,\psi\right\rangle
+\left\Vert \psi\right\Vert ^{2}\right\vert }-1\right\vert =\left\vert
\frac{\left\langle \psi,\psi\right\rangle +\left\Vert \psi\right\Vert ^{2}%
}{\left\vert \left\langle \psi,\psi\right\rangle +\left\Vert \psi\right\Vert
^{2}\right\vert }-\frac{\left\langle \varphi,\varphi\right\rangle +\left\Vert
\varphi\right\Vert ^{2}}{\left\vert \left\langle \varphi,\varphi\right\rangle
+\left\Vert \varphi\right\Vert ^{2}\right\vert }\right\vert \\
&  \qquad\leq\frac{\left\vert \left\langle \psi,\psi\right\rangle
-\left\langle \varphi,\varphi\right\rangle +\left\Vert \psi\right\Vert
^{2}-\left\Vert \varphi\right\Vert ^{2}\right\vert }{\left\vert \left\langle
\psi,\psi\right\rangle +\left\Vert \psi\right\Vert ^{2}\right\vert }\\
&  \qquad\quad+\frac{\left\vert \left\langle \varphi,\varphi\right\rangle
+\left\Vert \varphi\right\Vert ^{2}\right\vert }{\left\vert \left\langle
\psi,\psi\right\rangle +\left\Vert \psi\right\Vert ^{2}\right\vert \left\vert
\left\langle \varphi,\varphi\right\rangle +\left\Vert \varphi\right\Vert
^{2}\right\vert }\left\vert \left\vert \left\langle \varphi,\varphi
\right\rangle +\left\Vert \varphi\right\Vert ^{2}\right\vert -\left\vert
\left\langle \psi,\psi\right\rangle +\left\Vert \psi\right\Vert ^{2}%
\right\vert \right\vert .
\end{align*}
The inequalities%
\begin{align*}
\left\vert \left\langle \psi,\psi\right\rangle \right\vert  &  \leq\left\Vert
\psi\right\Vert ^{2}\leq\left(  \lambda+\beta\right)  ^{2}\left\Vert
p\right\Vert ^{2},\\
\left\Vert \varphi\right\Vert ^{2}  &  \leq\lambda^{2}\left\Vert p\right\Vert
^{2},\\
\left\vert \left\langle \varphi,\varphi\right\rangle +\left\Vert
\varphi\right\Vert ^{2}\right\vert  &  =2\left\Vert \varphi\right\Vert
^{2}\geq2\lambda^{2}\left\Vert p\right\Vert ^{2},\\
\left\Vert \psi\right\Vert ^{2}-\left\Vert \varphi\right\Vert ^{2}  &
\leq\left(  \left\Vert \psi\right\Vert +\left\Vert \varphi\right\Vert \right)
\left(  \left\Vert \psi\right\Vert -\left\Vert \varphi\right\Vert \right)
\leq\left\Vert p\right\Vert ^{2}\left(  2\lambda+\beta\right)  \beta,\\
\left\vert \left\langle \psi,\psi\right\rangle -\left\langle \varphi
,\varphi\right\rangle \right\vert  &  \leq\left(  \left\Vert \psi\right\Vert
+\left\Vert \varphi\right\Vert \right)  \left\Vert \psi-\varphi\right\Vert
\leq\left\Vert p\right\Vert ^{2}\left(  2\lambda+\beta\right)  \beta
\end{align*}
are derived by the reasoning: the first one follows by the same arguments as
in the proof of Lemma 3.6 in \cite{BoermMelenk}, the second one from the
definition of $\zeta$, the third one from the second one, the last two
inequalities from (\ref{normpsimnormvphi}). This leads to%
\begin{align*}
\left\vert \left\langle \psi,\psi\right\rangle +\left\Vert \psi\right\Vert
^{2}\right\vert  &  =\left\vert \left\langle \varphi,\varphi\right\rangle
+\left\Vert \varphi\right\Vert ^{2}\right\vert -\left\vert \left\langle
\psi,\psi\right\rangle +\left\Vert \psi\right\Vert ^{2}-\left\langle
\varphi,\varphi\right\rangle +\left\Vert \varphi\right\Vert ^{2}\right\vert \\
&  \geq2\left\Vert p\right\Vert ^{2}\left(  \lambda^{2}-\left(  2\lambda
+\beta\right)  \beta\right)  .
\end{align*}
The combination of these estimates leads to%
\[
\left\vert n_{\operatorname*{dp}}\left(  z\right)  -n_{\operatorname*{dp}%
}\left(  x_{z}\right)  \right\vert \leq\beta\left\Vert p\right\Vert \left(
1+2\frac{\left(  2\lambda+\beta\right)  \lambda}{\left(  \lambda^{2}-\left(
2\lambda+\beta\right)  \beta\right)  }\right)
\]
Now we use $\beta\leq\left(  \sqrt{\frac{3}{2}}-1\right)  \lambda$ to obtain%
\[
\left\vert n_{\operatorname*{dp}}\left(  z\right)  -n_{\operatorname*{dp}%
}\left(  x_{z}\right)  \right\vert \leq\left(  2\sqrt{6}+5\right)
\beta\left\Vert p\right\Vert .
\]
The estimate $\left\Vert p\right\Vert \leq\frac{1}{2}\max\left\{
\operatorname*{diam}\left(  B_{t}\right)  ,\operatorname*{diam}\left(
B_{s}\right)  \right\}  \leq\frac{\eta_{2}}{2}\operatorname*{dist}\left(
B_{t},B_{s}\right)  $ leads to%
\[
\operatorname{Re}\mathfrak{n}_{\operatorname*{dp}}\left(  z\right)
\geq\left\Vert d-x_{z}p\right\Vert -\left\vert n_{\operatorname*{dp}}\left(
z\right)  -n_{\operatorname*{dp}}\left(  x_{z}\right)  \right\vert \geq\left(
\lambda-\frac{\eta_{2}}{2}\left(  2\sqrt{6}+5\right)  \beta\right)  \left\Vert
p\right\Vert .
\]
Since $\beta\leq\frac{\lambda\left(  1-\eta_{3}\right)  }{\eta_{2}\left(
2\sqrt{6}+5\right)  }$ we have proved that%
\[
\operatorname{Re}\mathfrak{n}_{\operatorname*{dp}}\left(  z\right)  \geq
\frac{\eta_{3}+1}{2}\lambda\left\Vert p\right\Vert \geq\frac{1+\eta_{3}}%
{2}\operatorname*{dist}\left(  B_{t},B_{s}\right)
\]
holds. The combination of the estimates for $w\left(  \beta\right)  $,
$\chi_{\pm}\left(  \beta\right)  $ leads to%
\[
\sup_{z\in U_{\beta}}\left\vert G_{\operatorname*{dp}}\left(  z\right)
\right\vert \leq\frac{\exp\left(  \eta_{1}-\frac{1-\eta_{3}}{2}\left(
\operatorname{Re}\zeta\right)  \operatorname*{dist}\left(  B_{t},B_{s}\right)
\right)  }{\pi\operatorname*{dist}\left(  B_{t},B_{s}\right)  }%
\]

Note that $U_{\beta}$ contains the \textquotedblleft Bernstein
ellipse\textquotedblright:%
\[
\overline{\mathcal{D}}_{\rho}:=\left\{  z=x+\operatorname*{i}y:x,y\in
\mathbb{R},\left(  \frac{2x}{\rho+1/\rho}\right)  ^{2}+\left(  \frac{2y}%
{\rho-1/\rho}\right)  ^{2}\right\}  <1
\]
for%
\[
\rho:=\sqrt{\beta^{2}+1}+\beta\geq\sqrt{\widehat{\beta}^{2}+1}+\widehat{\beta
}>\widehat{\beta}+1=\rho_{0}.
\]
We have $\rho\geq\alpha\rho_{0}$ and $\alpha>1$ (cf. (\ref{defalpha})). We
know, e.g., from \cite[Lemma 3.11]{BoermMelenk}, that there exists
$q\in\mathbb{P}_{m}$ such that%
\[
\left\Vert G_{\operatorname*{dp}}-q\right\Vert _{\infty,\left[  -1,1\right]
}\leq\frac{2}{\rho_{0}-1}\rho^{-m}\frac{\exp\left(  \eta_{1}-\frac{1-\eta_{3}%
}{2}\left(  \operatorname{Re}\zeta\right)  \operatorname*{dist}\left(
B_{t},B_{s}\right)  \right)  }{\pi\operatorname*{dist}\left(  B_{t}%
,B_{s}\right)  }.
\]
By employing the Lebesgue constant $\Lambda_{m}$ for the \v{C}eby\v{s}ev
interpolation we conclude that%
\[
\left\Vert G_{\operatorname*{dp}}-\mathfrak{I}\left(  G_{\operatorname*{dp}%
}\right)  \right\Vert _{\infty,\left[  -1,1\right]  }\leq C_{1}\frac
{\alpha^{-m/2}\operatorname*{e}^{-\sigma\left(  \operatorname{Re}\zeta\right)
\operatorname*{dist}\left(  B_{t},B_{s}\right)  }}{4\pi\operatorname*{dist}%
\left(  B_{t},B_{s}\right)  }\rho_{0}^{-m}.
\]
The same arguments as in the proof of \cite[Corollary 3.14]{BoermMelenk} now
finishes the proof.%
%TCIMACRO{\TeXButton{End Proof}{\endproof}}%
%BeginExpansion
\endproof
%EndExpansion

By using the local error estimate we obtain the following consistency estimate.

\begin{theorem}
Let (\ref{sr}), (\ref{quconst}), (\ref{cpconst}) be satisfied. Let
$c\in\mathbb{R}^{3}$ and let the block $b=\left(  t,s\right)  $ satisfy the $%
%TCIMACRO{\TeXButton{boldeta}{\mbox{\boldmath$ \eta$}}}%
%BeginExpansion
\mbox{\boldmath$ \eta$}%
%EndExpansion
$-admissibility conditions (\ref{adm}) for some $0<\eta_{3}<1$. Then, there
exist constants $C_{\operatorname*{cons}}>0$ and $\sigma_{2}>0$ such that the
consistency estimate for the approximate sesquilinear form $a_{\zeta
}^{\mathcal{DH}^{2}}$ (cf. (\ref{h2approxazeta})) holds%
\[
\left\vert a_{\zeta}\left(  \varphi,\psi\right)  -a_{\zeta}^{\mathcal{DH}^{2}%
}\left(  \varphi,\psi\right)  \right\vert \leq C_{\operatorname*{cons}%
}C_{\Gamma}\frac{\rho_{0}^{-m}\operatorname*{e}^{-\sigma_{2}\left(
\operatorname{Re}\zeta\right)  h_{\mathcal{G}}}}{h_{\mathcal{G}}^{2}%
}\left\Vert \varphi\right\Vert _{L^{2}\left(  \Gamma\right)  }\left\Vert
\psi\right\Vert _{L^{2}\left(  \Gamma\right)  }.
\]
with%
\[
C_{\Gamma}:=\int_{\Gamma}\int_{\Gamma}\frac{1}{\left\Vert x-y\right\Vert
}d\Gamma_{y}d\Gamma_{x}.
\]

\end{theorem}

%

%TCIMACRO{\TeXButton{Proof}{\proof}}%
%BeginExpansion
\proof
%EndExpansion
Let $\varphi,\psi\in S$ with coefficient vectors $\left(  \varphi_{j}\right)
_{j=1}^{n}$, $\left(  \psi_{j}\right)  _{j=1}^{n}$ in their basis
representations (cf. (\ref{defpsiphi})). Then the difference $e_{\zeta
}:=a_{\zeta}-a_{\zeta}^{\mathcal{DH}^{2}}$ satisfies%
\[
\left\vert e_{\zeta}\left(  \varphi,\psi\right)  \right\vert \leq
\sum_{b=\left(  t,s\right)  \in\mathcal{P}_{\operatorname*{far}}}\left\vert
\sum_{i\in\hat{t}}\sum_{j\in\hat{s}}\overline{\psi_{i}}\varphi_{j}\int
_{\omega_{t}}\int_{\omega_{s}}E\left(  \zeta,x-y\right)  b_{i}\left(
x\right)  b_{j}\left(  y\right)  d\Gamma_{y}d\Gamma_{x}\right\vert
\]
for%
\[
E\left(  \zeta,x-y\right)  :=G\left(  \zeta,x-y\right)  -\sum_{%
%TCIMACRO{\TeXButton{boldmu}{\mbox{\boldmath$ \mu$}}}%
%BeginExpansion
\mbox{\boldmath$ \mu$}%
%EndExpansion
,%
%TCIMACRO{\TeXButton{boldnu}{\mbox{\boldmath$ \nu$}}}%
%BeginExpansion
\mbox{\boldmath$ \nu$}%
%EndExpansion
\in\mathbb{N}_{t}}\gamma_{%
%TCIMACRO{\TeXButton{boldmu}{\mbox{\boldmath$ \mu$}}}%
%BeginExpansion
\mbox{\boldmath$ \mu$}%
%EndExpansion
,%
%TCIMACRO{\TeXButton{boldnu}{\mbox{\boldmath$ \nu$}}}%
%BeginExpansion
\mbox{\boldmath$ \nu$}%
%EndExpansion
,c}^{b}\left(  \zeta\right)  \Phi_{%
%TCIMACRO{\TeXButton{boldmu}{\mbox{\boldmath$ \mu$}}}%
%BeginExpansion
\mbox{\boldmath$ \mu$}%
%EndExpansion
,c}^{t}\left(  \zeta,x\right)  \overline{\Phi_{%
%TCIMACRO{\TeXButton{boldmu}{\mbox{\boldmath$ \mu$}}}%
%BeginExpansion
\mbox{\boldmath$ \mu$}%
%EndExpansion
,c}^{s}}\left(  \zeta,y\right)  .
\]
Using the local error estimate (Thm. \ref{TheoErrLoc}) we obtain%
\[
\left\vert e_{\zeta}\left(  \varphi,\psi\right)  \right\vert \leq\frac
{C_{0}\rho_{0}^{-m}}{4\pi}\left\Vert \varphi\right\Vert _{L^{\infty}\left(
\Gamma\right)  }\left\Vert \psi\right\Vert _{L^{\infty}\left(  \Gamma\right)
}\sum_{b=\left(  t,s\right)  \in\mathcal{P}_{\operatorname*{far}}}\int
_{\omega_{t}}\int_{\omega_{s}}\frac{\operatorname*{e}^{-\sigma\left(
\operatorname{Re}\zeta\right)  \operatorname*{dist}\left(  B_{t},B_{s}\right)
}}{\operatorname*{dist}\left(  B_{t},B_{s}\right)  }d\Gamma_{y}d\Gamma_{x}.
\]
For all $\left(  x,y\right)  \in\omega_{t}\times\omega_{s}$, the standard
admissibility condition (cf. (\ref{adm_sing})) implies%
\[
\operatorname*{dist}\left(  B_{t},B_{s}\right)  \geq\left\Vert x-y\right\Vert
-\operatorname*{diam}B_{t}-\operatorname*{diam}B_{s}\geq\left\Vert
x-y\right\Vert -2\eta_{2}\operatorname*{dist}\left(  B_{t},B_{s}\right)
\]
and, in turn,%
\[
\operatorname*{dist}\left(  B_{t},B_{s}\right)  \geq\frac{1}{1+2\eta_{2}%
}\left\Vert x-y\right\Vert .
\]
We employ Remark \ref{Remdeltamin} to get
\[
\min\left\{  \operatorname*{dist}\left(  B_{t},B_{s}\right)  :\left(
t,s\right)  \text{ is admissible}\right\}  \geq\frac{h_{\min}}{\eta_{2}%
}\overset{\text{(\ref{quconst})}}{\geq}\frac{h_{\mathcal{G}}}%
{C_{\operatorname*{qu}}\eta_{2}}.
\]
Thus, for $\sigma_{1}:=\sigma/\left(  1+2\eta_{2}\right)  $ it holds%
\[
\left\vert e_{\zeta}\left(  \varphi,\psi\right)  \right\vert \leq\frac
{C_{0}\left(  1+2\eta_{2}\right)  \rho_{0}^{-m}}{4\pi}\left\Vert
\varphi\right\Vert _{L^{\infty}\left(  \Gamma\right)  }\left\Vert
\psi\right\Vert _{L^{\infty}\left(  \Gamma\right)  }\underset{\left\Vert
x-y\right\Vert \geq\frac{h_{\mathcal{G}}}{C_{\operatorname*{qu}}\eta_{2}}%
}{\int_{\Gamma}\int_{\Gamma}}\frac{\operatorname*{e}^{-\sigma_{1}\left(
\operatorname{Re}\zeta\right)  \left\Vert x-y\right\Vert }}{\left\Vert
x-y\right\Vert }d\Gamma_{y}d\Gamma_{x}%
\]
so that for $\sigma_{2}:=\sigma_{1}/\left(  C_{\operatorname*{qu}}\eta
_{2}\right)  $%
\[
\left\vert e_{\zeta}\left(  \varphi,\psi\right)  \right\vert \leq\frac
{C_{0}\left(  1+2\eta_{2}\right)  \rho_{0}^{-m}\operatorname*{e}^{-\sigma
_{2}\left(  \operatorname{Re}\zeta\right)  h_{\mathcal{G}}}}{4\pi}C_{\Gamma
}\left\Vert \varphi\right\Vert _{L^{\infty}\left(  \Gamma\right)  }\left\Vert
\psi\right\Vert _{L^{\infty}\left(  \Gamma\right)  }.
\]
The shape regularity and quasi-uniformity of the mesh implies (cf.
\cite[\S 4.4]{SauterSchwab2010}) that there exists a constant
$C_{\operatorname*{inv}}$ such that%
\[
\left\Vert \varphi\right\Vert _{L^{\infty}\left(  \Gamma\right)  }\leq
C_{\operatorname*{inv}}h_{\mathcal{G}}^{-1}\left\Vert \varphi\right\Vert
_{L^{2}\left(  \Gamma\right)  }%
\]
so that the assertion follows for $C_{\operatorname*{cons}}:=\frac
{C_{0}\left(  1+2\eta_{2}\right)  }{4\pi}C_{\Gamma}C_{\operatorname*{inv}}%
^{2}$.%
%TCIMACRO{\TeXButton{End Proof}{\endproof}}%
%BeginExpansion
\endproof
%EndExpansion

\section{Complexity\label{SecComplexity}}

In this section, we will estimate the complexity of the fast directional
$\mathcal{H}^{2}$-matrix approach for acoustic boundary integral operators
with complex frequency. In \cite{Boerm2015}, the complexity was analyzed for
purely complex frequencies $\zeta\in\operatorname*{i}\mathbb{R}$. Here we
generalize this theory to general complex frequencies $\zeta$ by taking into
account the modified admissibility condition \eqref{adm_paranew}. We will
derive explicit complexity estimates with respect to $\operatorname{Re}\zeta$
and $\operatorname{Im}\zeta$.

\subsection{Storage Requirements\label{SubSecStorage}}

\begin{remark}
[tridiagonal case]\label{RemTridiagCase}Let the boundary element mesh be
quasi-uniform and shape regular and the constants $c_{0}$, $\tilde{\sigma}$ as
in \ (\ref{defmbtilde}). Then the condition%
\[
\operatorname{Re}\zeta>c_{1}\sqrt{n}\log\frac{1}{\varepsilon}\quad\text{for
}c_{1}:=\frac{\tilde{C}c_{0}\eta_{2}}{\tilde{\sigma}}\text{\quad and\quad
}\tilde{C}:=\sqrt{\frac{C_{\operatorname*{sr}}}{C_{p}\left\vert \Gamma
\right\vert }}C_{\operatorname*{qu}}%
\]
implies that $m\left(  b\right)  =-1$ for all blocks $b\in\mathcal{P}%
_{\operatorname*{far}}$. As a consequence, the boundary element matrix
$\mathbf{K}\left(  \zeta\right)  $ can be replaced by its part, where the
kernel function is singular, i.e., $\mathbf{K}\left(  \zeta\right)
\approx\mathbf{K}^{0}\left(  \zeta\right)  =\left(  K_{i,j}^{0}\left(
\zeta\right)  \right)  _{i,j=1}^{n}$ with%
\[
K_{i,j}^{0}\left(  \zeta\right)  :=\left\{
\begin{array}
[c]{ll}%
K_{i,j}\left(  \zeta\right)  & \text{if }\operatorname*{dist}\left(
\omega_{i},\omega_{j}\right)  =0,\\
0 & \text{otherwise}%
\end{array}
\right.
\]
and $\omega_{i}=\operatorname*{supp}b_{i}$.
\end{remark}

%

%TCIMACRO{\TeXButton{Proof}{\proof}}%
%BeginExpansion
\proof
%EndExpansion
From Remark \ref{Remdeltamin} we obtain $\operatorname*{dist}\left(
B_{t},B_{s}\right)  \geq\frac{h_{\min}}{\eta_{2}}$ so that the following
implication holds (cf. (\ref{defmbtilde}))%
\begin{equation}
\frac{c_{0}\eta_{2}}{\tilde{\sigma}}\log\frac{1}{\varepsilon}<\left(
\operatorname{Re}\zeta\right)  h_{\min}\implies\tilde{m}_{b}=-1
\label{zerocond2}%
\end{equation}
For a quasi-uniform and shape regular boundary element mesh it holds%
\begin{equation}
\left\vert \Gamma\right\vert =\sum_{\tau\in\mathcal{G}}\left\vert
\tau\right\vert \leq C_{\operatorname*{sr}}\sum_{\tau\in\mathcal{G}}h_{\tau
}^{2}\leq C_{\operatorname*{sr}}h_{\mathcal{G}}^{2}\sharp\mathcal{G\leq}%
\frac{C_{\operatorname*{sr}}C_{\operatorname*{qu}}^{2}}{C_{p}}h_{\min}^{2}n
\label{nhest}%
\end{equation}
so that $h_{\min}\geq\tilde{C}n^{-1/2}$. The combination with (\ref{zerocond2}%
) leads to $\tilde{m}_{b}=-1$ for all $b\in\mathcal{P}_{\operatorname*{far}}$.
The definition of $m_{b}$ (cf. (\ref{defmb})) then finishes the proof.%
%TCIMACRO{\TeXButton{End Proof}{\endproof}}%
%BeginExpansion
\endproof
%EndExpansion

\begin{remark}
[sectorial case]\label{Remsect}For $\zeta\in\mathbb{C}_{>0}$ with $\left\vert
\operatorname{Im}\zeta\right\vert \leq\alpha\operatorname{Re}\zeta$ and some
$\alpha>0$ the condition%
\[
\max\left\{  \operatorname*{diam}\left(  B_{t}\right)  ,\operatorname*{diam}%
\left(  B_{s}\right)  \right\}  \leq\sqrt{\frac{\eta_{3}}{\alpha}%
}\operatorname*{dist}\left(  B_{t},B_{s}\right)
\]
is stronger than the condition (\ref{adm_paranew}) and can be absorbed into
the condition (\ref{adm_sing}) by adjusting $\eta_{2}\leftarrow\tilde{\eta
}_{2}:=\min\left\{  \eta_{2},\sqrt{\frac{\eta_{3}}{\alpha}}\right\}  $. Since
condition (\ref{adm_angle}) is understood as a condition on the number and
choice of directions in $\mathcal{D}_{\ell}$, the number of elements in the
minimal partition $\mathcal{P}$ of $\mathcal{I}\times\mathcal{I}$ is bounded
by the number of elements in a partition $\mathcal{\tilde{P}}$ where
conditions (\ref{adm_sing}) and (\ref{adm_paranew}) are replaced by the
condition $\max\left\{  \operatorname*{diam}\left(  B_{t}\right)
,\operatorname*{diam}\left(  B_{s}\right)  \right\}  \leq\tilde{\eta}%
_{2}\operatorname*{dist}\left(  B_{t},B_{s}\right)  $. This is the standard
admissibility condition for the Laplacian and estimates of the form
$\sharp\mathcal{\tilde{P}\leq C}n$ are well known (see, e.g., \cite{Sauter00},
\cite{SauterSchwab2010}).
\end{remark}

Next, we will estimate the number of elements in the minimal partition
$\mathcal{P}$ for the case $\zeta\in\mathbb{C}_{>0}$ with%
\[
\left\vert \operatorname{Im}\zeta\right\vert >\alpha\operatorname{Re}%
\zeta\quad\text{for }\alpha\text{ as in Rem. \ref{Remsect}.}%
\]
Note that this condition implies that $\left\vert \zeta\right\vert $ and
$\left\vert \operatorname{Im}\zeta\right\vert $ are \textquotedblleft
equivalent\textquotedblright:%
\[
\left\vert \operatorname{Im}\zeta\right\vert \leq\left\vert \zeta\right\vert
\leq\sqrt{1+\alpha^{-2}}\left\vert \operatorname{Im}\zeta\right\vert .
\]

The theory in this Section is a generalization of the one in \cite[Section
5]{Boerm2015}, adapted to our new admissibility condition \eqref{adm_paranew}
and our goal is to derive estimates which are explicit in all relevant
parameters, in particular, with respect to $\operatorname{Im}\zeta$,
$\operatorname{Re}\zeta$, $n$, and certain geometric parameters which we will
introduce next.

As in \cite{Boerm2015} we assume that there exist a reference box $B_{\ell}$
for each level $\ell$ and constants $\rho_{\operatorname*{ref}}>1$,
$C_{\operatorname*{sb}},C_{\operatorname*{sn}}\geq1$, $c_{\operatorname*{ref}%
},C_{\operatorname*{bp}},C_{\operatorname*{bb}},C_{\operatorname*{rs}%
},C_{\operatorname*{ov}},C_{\operatorname*{un}}>0$ such that:
\begin{align}
&  \exists d_{t}\in\mathbb{R}^{3}\text{\quad s.t.\quad}B_{t}=B_{\ell}%
+d_{t},\quad\mbox{ for all clusters }t\in\mathcal{T}_{\mathcal{I}}^{(\ell
)},\label{10}\\
&  \operatorname*{diam}(B_{t})\leq C_{\operatorname*{sb}}\operatorname*{diam}%
(B_{t^{\prime}}),\quad\mbox{for all }t\in\mathcal{T}_{\mathcal{I}},t^{\prime
}\in\operatorname*{sons}\left(  t\right)  ,\label{11}\\
&  \mbox{\#}\operatorname*{sons}(t)\leq C_{\operatorname*{sn}},\quad
\mbox{\#}\operatorname*{sons}(t)\neq1,\quad\mbox{for all }t\in\mathcal{T}%
_{\mathcal{I}},\label{12}\\
&  \left\vert \Gamma\cap\mathcal{B}(x,r)\right\vert \leq C_{\operatorname*{bp}%
}r^{2},\quad\mbox{ for all }x\in\mathbb{R}^{3},r\geq0,\label{13a}\\
&  \operatorname*{diam}\nolimits^{2}(B_{t})\leq C_{\operatorname*{bb}%
}\left\vert B_{t}\cap\Gamma\right\vert ,\quad\mbox{ for all }t\in
\mathcal{T}_{\mathcal{I}},\label{13b}\\
&  \mbox{\#}\{t\in\mathcal{T}_{\mathcal{I}}^{(\ell)}\ :\ x\in B_{t}\}\leq
C_{\operatorname*{ov}},\quad\mbox{ for all }x\in\Omega,\ell\in\mathbb{N}%
_{0},\label{14}\\
&  C_{\operatorname*{rs}}^{-1}\left(  k_{L}+1\right)  \leq\mbox{\#}\hat{t}\leq
C_{\operatorname*{rs}}\left(  k_{L}+1\right)  ,\quad
\mbox{ for all leaves }t\in\mathcal{L}_{\mathcal{I}}\cap\mathcal{T}_{\ell
},\label{15b}\\
&  c_{\operatorname*{ref}}\rho_{\operatorname*{ref}}^{L-\ell}\delta_{L}%
\leq\delta_{\ell}\quad\forall0\leq\ell\leq L:=\operatorname*{depth}%
\mathcal{T}_{\mathcal{I}}\text{ (with }\delta_{\ell}\text{ as in
(\ref{maxcldiam})).} \label{delta_ell_cond}%
\end{align}

\begin{lemma}
[Sparsity]\label{lemma:inadmin} Let (\ref{10}-\ref{14}) hold. For every
cluster $t\in\mathcal{T}_{\mathcal{I}}$, the sets $\mathcal{P}%
_{\operatorname*{left}}\left(  t\right)  $, $\mathcal{P}%
_{\operatorname*{right}}\left(  t\right)  $ as in (\ref{leftright}) satisfies%
\begin{equation}
\max\left\{  \sharp\mathcal{P}_{\operatorname*{left}}\left(  t\right)
,\sharp\mathcal{P}_{\operatorname*{right}}\left(  t\right)  \right\}  \leq
\hat{C}_{\operatorname*{sp}}R_{t}^{2} \label{ninadm}%
\end{equation}
with $\hat{C}_{\operatorname*{sp}}:=C_{\operatorname*{sn}}%
C_{\operatorname*{sb}}^{2}C_{\operatorname*{bb}}C_{\operatorname*{ov}%
}C_{\operatorname*{bp}}$ and%
\begin{equation}
R_{t}:=\frac{3}{2}+\max\left\{  \frac{1}{\eta_{2}},r_{t}\right\}
\quad\text{and\quad}r_{t}:=\min\left\{  \frac{\left\vert \operatorname{Im}%
\zeta\right\vert }{\eta_{2}}\operatorname*{diam}B_{t},\sqrt{\frac{\left\vert
\operatorname{Im}\zeta\right\vert }{\eta_{3}\operatorname{Re}\zeta}}\right\}
. \label{DefRtrt}%
\end{equation}

\end{lemma}

%

%TCIMACRO{\TeXButton{Proof}{\proof}}%
%BeginExpansion
\proof
%EndExpansion
We prove the estimate only for $\sharp\mathcal{P}_{\operatorname*{right}%
}\left(  t\right)  $ while the proof for $\sharp\mathcal{P}%
_{\operatorname*{left}}\left(  t\right)  $ follows verbatim.

Let $t\in\mathcal{T}_{I}$ and $s\in\mathcal{P}_{\operatorname*{right}%
}^{\operatorname*{near}}\left(  t\right)  $ (cf. (\ref{leftright})). Since $t$
and $s$ belong to the same tree level we have $\max\left\{
\operatorname*{diam}B_{t},\operatorname*{diam}B_{s}\right\}
=\operatorname*{diam}B_{t}$. Then, for any $z\in B_{s}$ the estimate%
\begin{equation}
\left\Vert z-M_{t}\right\Vert \leq\operatorname*{diam}B_{s}%
+\operatorname*{dist}\left(  B_{t},B_{s}\right)  +\frac{1}{2}%
\operatorname*{diam}B_{t} \label{distsplit}%
\end{equation}
holds. Since the block $\left(  t,s\right)  $ is non-admissible one of the
conditions (\ref{adm_sing}), (\ref{adm_paranew}) must be violated.

\textbf{Case 1.} Let condition (\ref{adm_sing}) be violated. Then,%
\[
\left\Vert z-M_{t}\right\Vert <\left(  \frac{3}{2}+\frac{1}{\eta_{2}}\right)
\operatorname*{diam}B_{t}.
\]

\textbf{Case 2.} Let condition (\ref{adm_paranew}) be violated but condition
(\ref{adm_sing}) be valid. We set $\lambda_{t,s}:=\max\left\{  1,\frac
{\eta_{3}}{\eta_{2}}\left(  \operatorname{Re}\zeta\right)
\operatorname*{dist}\left(  B_{t},B_{s}\right)  \right\}  $ and obtain by
combining (\ref{adm_sing}) with the negation of (\ref{adm_paranew})%
\begin{equation}
\left\vert \operatorname{Im}\zeta\right\vert \left(  \operatorname*{diam}%
B_{t}\right)  ^{2}>\lambda_{t,s}\eta_{2}\operatorname*{dist}\left(
B_{t},B_{s}\right)  \geq\lambda_{t,s}\operatorname*{diam}B_{t} \label{imdiam2}%
\end{equation}
so that\quad$\left\vert \operatorname{Im}\zeta\right\vert \operatorname*{diam}%
B_{t}>\lambda_{t,s}$. The left inequality in (\ref{imdiam2}) can be split into%
\begin{align*}
\operatorname*{dist}\left(  B_{t},B_{s}\right)   &  <\frac{\left\vert
\operatorname{Im}\zeta\right\vert }{\eta_{2}}\left(  \operatorname*{diam}%
B_{t}\right)  ^{2}\\
\operatorname*{dist}\left(  B_{t},B_{s}\right)   &  <\sqrt{\frac{\left\vert
\operatorname*{Im}\zeta\right\vert }{\eta_{3}\operatorname{Re}\zeta}%
}\operatorname*{diam}B_{t}%
\end{align*}
so that%
\[
\operatorname*{dist}\left(  B_{t},B_{s}\right)  <r_{t}\operatorname*{diam}%
B_{t}%
\]
for $r_{t}$ as in (\ref{DefRtrt}). The combination of this with
(\ref{distsplit}) leads to%
\[
\left\Vert z-M_{t}\right\Vert \leq\left(  r_{t}+\frac{3}{2}\right)
\operatorname*{diam}B_{t}.
\]
The distance estimates in Case 1 and Case 2 lead to%
\[
\left\Vert z-M_{t}\right\Vert <\min\left\{  \left(  \frac{3}{2}+\frac{1}%
{\eta_{2}}\right)  \operatorname*{diam}B_{t},\left(  r_{t}+\frac{3}{2}\right)
\operatorname*{diam}B_{t}\right\}  \quad\forall z\in B_{s}%
\]
for non-admissible pairs of clusters and allow for an estimate of the
cardinality. The cluster $s\in\mathcal{P}_{\operatorname*{right}%
}^{\operatorname*{near}}\left(  t\right)  $ is contained in a ball with center
$M_{t}$ and radius $R_{t}$ as in (\ref{DefRtrt}). By the same arguments as in
the proof of \cite[Lem. 2]{Boerm2015} we obtain%
\[
\sharp\mathcal{P}_{\operatorname*{right}}^{\operatorname*{near}}\left(
t\right)  \leq C_{\operatorname*{bb}}C_{\operatorname*{ov}}%
C_{\operatorname*{bp}}R_{t}^{2}.
\]

It remains to estimate the cardinality of $\sharp\mathcal{P}%
_{\operatorname*{right}}\left(  t\right)  $. If $t$ is the root of
$\mathcal{T}_{\mathcal{I}}$ we have $\mathcal{P}_{\operatorname*{right}%
}\left(  t\right)  =\mathcal{P}_{\operatorname*{left}}\left(  t\right)
=\left\{  t\right\}  $ and the cardinalities equal $1$. For $t\in
\mathcal{T}_{\mathcal{I}}$ being not the root we denote by $t^{+}$ the father
of $t$ which, by construction, is non-admissible. Hence,
\[
\sharp\mathcal{P}_{\operatorname*{right}}\left(  t\right)  \leq
C_{\operatorname*{sn}}\sharp\mathcal{P}_{\operatorname*{right}}%
^{\operatorname*{near}}\left(  t^{+}\right)  \leq C_{\operatorname*{sn}%
}C_{\operatorname*{bb}}C_{\operatorname*{ov}}C_{\operatorname*{bp}}R_{t^{+}%
}^{2}.
\]
The final estimate follows from%
\[
r_{t^{+}}\leq\min\left\{  C_{\operatorname*{sb}}\frac{\left\vert
\operatorname{Im}\zeta\right\vert }{\eta_{2}}\operatorname*{diam}B_{t}%
,\sqrt{\frac{\left\vert \operatorname{Im}\zeta\right\vert }{\eta
_{3}\operatorname{Re}\zeta}}\right\}  \leq C_{\operatorname*{sb}}r_{t}.
\]%
%TCIMACRO{\TeXButton{End Proof}{\endproof}}%
%BeginExpansion
\endproof
%EndExpansion

\begin{corollary}
[Nearfield sparsity with resolution condition]\label{LemSparsRes}For the
boundary element discretization we assume (\ref{sr}), (\ref{quconst}),
(\ref{cpconst}). Let (\ref{10}-\ref{14}) hold and assume the \emph{resolution
condition}%
\begin{equation}
\left\vert \operatorname{Im}\zeta\right\vert \operatorname*{diam}B_{t}\leq
C_{\operatorname*{res}}\quad\forall t\in\mathcal{L}_{\mathcal{I}%
}.\label{rescond}%
\end{equation}
Then, there exists a constant $C_{\operatorname*{near}}^{\sharp}$ such that
cardinality of the near field $\mathcal{P}_{\operatorname*{near}}$ as in
(\ref{defPnearPfar}) can be estimated%
\[
\sharp\mathcal{P}_{\operatorname*{near}}\leq C_{\operatorname*{near}}^{\sharp
}\frac{n}{k_{L}+1}.
\]

\end{corollary}

%

%TCIMACRO{\TeXButton{Proof}{\proof}}%
%BeginExpansion
\proof
%EndExpansion
Let $b=\left(  t,s\right)  \in\mathcal{P}_{\operatorname*{near}}$. Then, $t$
is a leaf or $s$ is a leaf. In the first case, it holds $s\in\mathcal{P}%
_{\operatorname*{right}}^{\operatorname*{near}}\left(  t\right)  $ and in the
second $t\in\mathcal{P}_{\operatorname*{left}}^{\operatorname*{near}}\left(
s\right)  $. Let $t$ be a leaf. Then we employ (\ref{rescond}) to obtain%
\begin{equation}
r_{t}\overset{\text{(\ref{DefRtrt})}}{\leq}\frac{C_{\operatorname*{res}}}%
{\eta_{2}}\quad\text{and\quad}R_{t}\leq\frac{3}{2}+\frac{\max\left\{
C_{\operatorname*{res}},1\right\}  }{\eta_{2}}. \label{estRcapt}%
\end{equation}
The combination with (\ref{ninadm}) leads to%
\[
\mathcal{P}_{\operatorname*{right}}^{\operatorname*{near}}\left(  t\right)
\leq\hat{C}_{\operatorname*{sp}}R_{t}^{2}.
\]
The proof of $\mathcal{P}_{\operatorname*{left}}^{\operatorname*{near}}\left(
s\right)  \leq\hat{C}_{\operatorname*{sp}}R_{s}^{2}$ in case that $s$ is a
leaf is verbatim Hence,%
\[
\max\left\{  \sharp\mathcal{P}_{\operatorname*{left}}^{\operatorname*{near}%
}\left(  t\right)  ,\sharp\mathcal{P}_{\operatorname*{right}}%
^{\operatorname*{near}}\left(  t\right)  \right\}  \leq\hat{C}%
_{\operatorname*{sp}}R_{t}^{2},\qquad\forall t\in\mathcal{T}_{\mathcal{I}}%
\]
and%
\[
\sharp\mathcal{P}_{\operatorname*{near}}\leq\sum_{t\in\mathcal{L}%
_{\mathcal{I}}}\left(  \sharp\mathcal{P}_{\operatorname*{left}}%
^{\operatorname*{near}}\left(  t\right)  +\sharp\mathcal{P}%
_{\operatorname*{right}}^{\operatorname*{near}}\left(  t\right)  \right)
\leq2\hat{C}_{\operatorname*{sp}}\sum_{t\in\mathcal{L}_{\mathcal{I}}}R_{t}%
^{2}\overset{\text{(\ref{estRcapt})}}{\leq}2\hat{C}_{\operatorname*{sp}%
}\left(  \frac{3}{2}+\frac{\max\left\{  C_{\operatorname*{res}},1\right\}
}{\eta_{2}}\right)  ^{2}\left(  \sharp\mathcal{L}_{\mathcal{I}}\right)
\]
For $\sharp\mathcal{L}_{\mathcal{I}}$, we obtain%
\begin{equation}
n=\sharp\mathcal{I=}\sum_{t\in\mathcal{L}_{\mathcal{I}}}\sharp\hat{t}%
\overset{\text{(\ref{15b})}}{\geq}C_{\operatorname*{rs}}^{-1}\left(
k_{L}+1\right)  \sharp\mathcal{L}_{\mathcal{I}} \label{estleaves}%
\end{equation}
and the assertion follows with%
\[
C_{\operatorname*{near}}^{\sharp}:=2\hat{C}_{\operatorname*{sp}}\left(
\frac{3}{2}+\frac{\max\left\{  C_{\operatorname*{res}},1\right\}  }{\eta_{2}%
}\right)  ^{2}C_{\operatorname*{rs}}.
\]%
%TCIMACRO{\TeXButton{End Proof}{\endproof}}%
%BeginExpansion
\endproof
%EndExpansion

The next estimate of the number of clusters is proven in \cite[Lem.
3]{Boerm2015} and carries over to our case, since it does not involve the
admissibility conditions \eqref{adm}.

\begin{lemma}
[Clusters]\label{lemma:clusters} Let \eqref{10} , \eqref{12}, \eqref{13b},
\eqref{14} and \eqref{15b} hold. Then%
\begin{equation}
\mbox{\#}\mathcal{T}_{\ell}\leq C_{\operatorname*{lv}}\frac{\left\vert
\Gamma\right\vert }{\operatorname*{diam}\nolimits^{2}(B_{\ell})}%
,\qquad\mbox{ for all }\ell\in\mathbb{N}_{0} \label{18a}%
\end{equation}
with $C_{\operatorname*{lv}}:=\max\{C_{\operatorname*{bb}}%
C_{\operatorname*{ov}},2C_{\operatorname*{rs}}\}$.
\end{lemma}

Next we estimate the cardinality of the cluster basis. For each $t\in
\mathcal{T}_{\ell}$ and $c\in\mathcal{D}_{\ell}$ we define%
\[
\mathcal{P}_{\operatorname*{right}}^{\operatorname*{far}}\left(  t,c\right)
:=\left\{  s\in\mathcal{T}_{\mathcal{I}}\mid b=\left(  t,s\right)
\in\mathcal{P}_{\operatorname*{adm}}\wedge c\left(  b\right)  =c\right\}
\]
and observe that%
\[%
%TCIMACRO{\dbigcup \limits_{c\in\mathcal{D}_{\ell}}}%
%BeginExpansion
{\displaystyle\bigcup\limits_{c\in\mathcal{D}_{\ell}}}
%EndExpansion
\mathcal{P}_{\operatorname*{right}}^{\operatorname*{far}}\left(  t,c\right)
=\mathcal{P}_{\operatorname*{right}}^{\operatorname*{far}}\left(  t\right)
\]
holds.

\begin{lemma}
[Block and cluster sums]\label{lemma:basis}Let the set of directions be
constructed according to \cite[Rem. 3]{borm2018hybrid}. Under assumptions
(\ref{cpconst}), (\ref{Cwidth}), (\ref{10}), (\ref{12}-\ref{15b}),
(\ref{delta_ell_cond}) and assuming that $L:=\operatorname*{depth}\left(
\mathcal{T}_{\mathcal{I}}\right)  >0$, there exists a constant $C_{\sharp}$
such that%
\begin{align}
&  \sum_{t\in\mathcal{T}_{\mathcal{I}}}\sum_{c\in\mathcal{D}%
_{\operatorname{level}\left(  t\right)  }}\sharp\mathcal{P}%
_{\operatorname*{right}}\left(  t,c\right) \label{lemma:basisa}\\
&  \qquad\leq C_{\sharp}\left(  \frac{n}{\eta_{2}^{2}\left(  k_{L}+1\right)
}+\min\left\{  \left(  1+L\right)  \left(  \frac{\left\vert \operatorname{Im}%
\zeta\right\vert }{\eta_{2}}\right)  ^{2},\frac{n}{\eta_{3}\left(
k_{L}+1\right)  }\frac{\left\vert \operatorname{Im}\zeta\right\vert
}{\operatorname{Re}\zeta}\right\}  \right)  .\nonumber
\end{align}
There exists $C_{\operatorname*{ffc}}>0$ such that the total number of basis
farfield coefficients (cf. (\ref{basisffc})) is bounded from above by%
\begin{equation}
C_{\operatorname*{ffc}}k_{L}\left(  n+k_{L}\left(  \operatorname{Im}%
\zeta\right)  ^{2}\right)  \label{lemma:basisb}%
\end{equation}
The total number of expansion coefficients (cf. (\ref{expcoeff})) is bounded
from above by%
\begin{equation}
C_{\sharp}k_{L}\left(  \frac{n}{\eta_{2}^{2}}+\min\left\{  k_{L}\left(
\left(  1+L\right)  \frac{\left\vert \operatorname{Im}\zeta\right\vert }%
{\eta_{2}}\right)  ^{2},\frac{n}{\eta_{3}}\frac{\left\vert \operatorname{Im}%
\zeta\right\vert }{\operatorname{Re}\zeta}\right\}  \right)  .
\label{lemma:basisc}%
\end{equation}

\end{lemma}

%

%TCIMACRO{\TeXButton{Proof}{\proof}}%
%BeginExpansion
\proof
%EndExpansion
\textbf{Part 1. Estimate of the total number of blocks}

We follow the arguments in the proof of \cite[Lem. 8]{borm2018hybrid}. Let
$L:=\operatorname*{depth}\left(  \mathcal{T}_{\mathcal{I}}\right)  $. The
combination of (\ref{18a}) with (\ref{ninadm}) leads to%
\[
\sum_{t\in\mathcal{T}_{\mathcal{I}}}\sum_{c\in\mathcal{D}%
_{\operatorname{level}\left(  t\right)  }}\sharp\mathcal{P}%
_{\operatorname*{right}}\left(  t,c\right)  \leq\hat{C}_{\operatorname*{sp}%
}\sum_{t\in\mathcal{T}_{\mathcal{I}}}R_{t}^{2}\leq\hat{C}_{\operatorname*{sp}%
}C_{\operatorname*{lv}}\left\vert \Gamma\right\vert \sum_{\ell=0}^{L}%
\frac{R_{\ell}^{2}}{\operatorname*{diam}\nolimits^{2}B_{\ell}}%
\]
with%
\begin{equation}
R_{\ell}^{2}\leq2\left(  \left(  \frac{3}{2}+\frac{1}{\eta_{2}}\right)
^{2}+r_{\ell}^{2}\right)  \label{estRl2}%
\end{equation}
and we obtain%
\begin{align*}
&  \sum_{t\in\mathcal{T}_{\mathcal{I}}}\sum_{c\in\mathcal{D}%
_{\operatorname{level}\left(  t\right)  }}\sharp\mathcal{P}%
_{\operatorname*{right}}\left(  t,c\right) \\
&  \quad\leq2\hat{C}_{\operatorname*{sp}}C_{\operatorname*{lv}}\left\vert
\Gamma\right\vert \left(  \sum_{\ell=0}^{L}\frac{\left(  \frac{3}{2}+\frac
{1}{\eta_{2}}\right)  ^{2}}{\operatorname*{diam}\nolimits^{2}B_{\ell}}%
+\sum_{\ell=0}^{L}\min\left\{  \left(  \frac{\left\vert \operatorname{Im}%
\zeta\right\vert }{\eta_{2}}\right)  ^{2},\frac{\left\vert \operatorname{Im}%
\zeta\right\vert }{\eta_{3}\left(  \operatorname{Re}\zeta\right)
\operatorname*{diam}^{2}B_{\ell}}\right\}  \right)  .
\end{align*}
We have for all $0\leq\ell\leq L$ and $t\in\mathcal{L}_{\mathcal{I}}%
\cap\mathcal{T}_{\ell}$ the estimate%
\begin{align}
\operatorname*{diam}\nolimits^{2}B_{\ell}  &  \overset{\text{(\ref{Cwidth})}%
}{\geq}c_{\operatorname*{vol}}\left\vert \omega_{t}\right\vert \overset
{\text{(\ref{sr})}}{\geq}c_{\operatorname*{sr}}\frac{\sharp\hat{t}%
}{C_{\operatorname*{loc},p}+1}h_{\min}^{2}\overset{\text{(\ref{quconst})}%
}{\geq}\frac{c_{\operatorname*{sr}}}{C_{\operatorname*{qu}}^{2}}\frac
{\sharp\hat{t}}{C_{\operatorname*{loc},p}+1}h_{\mathcal{G}}^{2}%
\label{2ndineququ}\\
&  \overset{\text{(\ref{nhest})}}{\geq}\frac{c_{\operatorname*{sr}}%
}{C_{\operatorname*{qu}}^{2}C_{\operatorname*{sr}}}\frac{\left\vert
\Gamma\right\vert }{C_{\operatorname*{loc},p}+1}\frac{\sharp\hat{t}}%
{\sharp\mathcal{G}}\overset{\text{(\ref{cpconst})}}{\geq}\frac
{c_{\operatorname*{sr}}}{C_{\operatorname*{qu}}^{2}C_{\operatorname*{sr}}%
}\frac{C_{p}\left\vert \Gamma\right\vert }{C_{\operatorname*{loc},p}+1}%
\frac{\sharp\hat{t}}{n}\overset{\text{(\ref{15b})}}{\geq}c_{\mathcal{L}%
}\left\vert \Gamma\right\vert \frac{k_{L}+1}{n}\nonumber
\end{align}
with%
\[
c_{\mathcal{L}}:=\frac{c_{\operatorname*{sr}}}{C_{\operatorname*{qu}}%
^{2}C_{\operatorname*{sr}}C_{\operatorname*{rs}}}\frac{C_{p}}%
{C_{\operatorname*{loc},p}+1}.
\]

We use (\ref{delta_ell_cond}) and get by a geometric sum argument%
\[
\sum_{\ell=0}^{L}\frac{1}{\operatorname*{diam}\nolimits^{2}B_{\ell}}\leq
c_{\operatorname*{ref}}^{2}\frac{1}{\rho_{\operatorname*{ref}}^{2}-1}\frac
{1}{\delta_{L}^{2}}\leq\frac{c_{B}}{\left\vert \Gamma\right\vert }\frac
{n}{k_{L}+1}\quad\text{with\quad}c_{B}:=\frac{c_{\operatorname*{ref}}^{2}%
}{c_{\mathcal{L}}\left(  \rho_{\operatorname*{ref}}^{2}-1\right)  }.
\]
We end up with the estimate%
\begin{align*}
&  \sum_{t\in\mathcal{T}_{\mathcal{I}}}\sum_{c\in\mathcal{D}%
_{\operatorname{level}\left(  t\right)  }}\sharp\mathcal{P}%
_{\operatorname*{right}}\left(  t,c\right) \\
&  \quad\leq2\hat{C}_{\operatorname*{sp}}C_{\operatorname*{lv}}\left(  \left(
\frac{3}{2}+\frac{1}{\eta_{2}}\right)  ^{2}c_{B}\frac{n}{k_{L}+1}+\min\left\{
\left\vert \Gamma\right\vert \left(  L+1\right)  \left(  \frac{\left\vert
\operatorname{Im}\zeta\right\vert }{\eta_{2}}\right)  ^{2},c_{B}%
\frac{\left\vert \operatorname{Im}\zeta\right\vert }{\eta_{3}\left(
\operatorname{Re}\zeta\right)  }\frac{n}{k_{L}+1}\right\}  \right)  .
\end{align*}
Hence, (\ref{lemma:basisa}) follows with%
\[
C_{\sharp}:=2\hat{C}_{\operatorname*{sp}}C_{\operatorname*{lv}}\max\left\{
\left(  \frac{3}{2}\eta_{2}+1\right)  ^{2}c_{B},\left\vert \Gamma\right\vert
\right\}  .
\]

\textbf{Part 2. Estimate of the total number of basis farfield coefficients}
$J_{%
%TCIMACRO{\TeXButton{boldmu}{\mbox{\boldmath$ \mu$}}}%
%BeginExpansion
\mbox{\boldmath$ \mu$}%
%EndExpansion
,c}^{t}\left(  \zeta,b_{j}\right)  $.

From \cite[(16)]{borm2018hybrid} we conclude that the construction of the set
of directions as in \cite[Rem. 3]{borm2018hybrid} implies%
\[
\sharp D_{\ell}\leq C_{\operatorname*{di}}\left(  1+\left(  \operatorname*{Im}%
\zeta\right)  ^{2}\operatorname*{diam}\nolimits^{2}B_{\ell}\right)
\quad\forall0\leq\ell\leq L:=\operatorname*{depth}\mathcal{T}_{\mathcal{I}}%
\]
for some constant $C_{\operatorname*{di}}>0$. Furthermore, we have%
\[
\sharp\left\{  j\in\hat{t}\right\}  \leq C_{\operatorname*{loc},p}\sharp
\hat{t}.
\]
The coefficients $J_{%
%TCIMACRO{\TeXButton{boldmu}{\mbox{\boldmath$ \mu$}}}%
%BeginExpansion
\mbox{\boldmath$ \mu$}%
%EndExpansion
,c}^{t}\left(  \zeta,b_{j}\right)  $ must be computed only for leaves
$t\in\mathcal{L}_{\mathcal{I}}$ and we obtain%
\begin{align*}
\sum_{t\in\mathcal{L}_{\mathcal{I}}}\sum_{%
%TCIMACRO{\TeXButton{boldmu}{\mbox{\boldmath$ \mu$}}}%
%BeginExpansion
\mbox{\boldmath$ \mu$}%
%EndExpansion
\in\mathbb{N}_{t}}\sum_{c\in\mathcal{D}_{t}}\sum_{j\in\hat{t}}1  &  \leq
C_{\operatorname*{loc},p}\sum_{t\in\mathcal{L}_{\mathcal{I}}}\sum_{%
%TCIMACRO{\TeXButton{boldmu}{\mbox{\boldmath$ \mu$}}}%
%BeginExpansion
\mbox{\boldmath$ \mu$}%
%EndExpansion
\in\mathbb{N}_{t}}\sum_{c\in\mathcal{D}_{\operatorname{level}\left(  t\right)
}}\sharp\hat{t}\\
&  \leq C_{\operatorname*{loc},p}C_{\operatorname*{di}}\sum_{t\in
\mathcal{L}_{\mathcal{I}}}k_{t}\left(  \sharp\hat{t}\right)  \left(  1+\left(
\operatorname*{Im}\zeta\right)  ^{2}\operatorname*{diam}\nolimits^{2}%
B_{\operatorname{level}\left(  t\right)  }\right) \\
&  \leq2C_{\operatorname*{loc},p}C_{\operatorname*{di}}k_{L}\left(
n+C_{\operatorname*{rs}}k_{L}\left(  \operatorname*{Im}\zeta\right)  ^{2}%
\sum_{t\in\mathcal{L}_{\mathcal{I}}}\operatorname*{diam}\nolimits^{2}%
B_{\operatorname{level}\left(  t\right)  }\right)  .
\end{align*}
It holds%
\[
\sum_{t\in\mathcal{L}_{\mathcal{I}}}\operatorname*{diam}\nolimits^{2}B_{\ell
}\overset{\text{(\ref{Cwidth})}}{\leq}C_{\operatorname*{vol}}\sum
_{t\in\mathcal{L}_{\mathcal{I}}}\left\vert \omega_{t}\right\vert \leq
C_{\operatorname*{vol}}C_{p}\left\vert \Gamma\right\vert .
\]
This allows to estimate%
\[
\sum_{t\in\mathcal{L}_{\mathcal{I}}}\sum_{%
%TCIMACRO{\TeXButton{boldmu}{\mbox{\boldmath$ \mu$}}}%
%BeginExpansion
\mbox{\boldmath$ \mu$}%
%EndExpansion
\in\mathbb{N}_{t}}\sum_{c\in\mathcal{D}_{\ell}}\sum_{j\in\hat{t}}%
1\leq2C_{\operatorname*{loc},p}C_{\operatorname*{di}}k_{L}\left(
n+C_{\operatorname*{rs}}C_{\operatorname*{vol}}C_{p}\left\vert \Gamma
\right\vert k_{L}\left(  \operatorname*{Im}\zeta\right)  ^{2}\right)
\]
from which (\ref{lemma:basisb}) follows with%
\[
C_{\operatorname*{ffc}}:=2C_{\operatorname*{loc},p}C_{\operatorname*{di}}%
\max\left\{  1,C_{\operatorname*{rs}}C_{\operatorname*{vol}}C_{p}\left\vert
\Gamma\right\vert \right\}  .
\]

\textbf{Part 3. Estimate the total number of expansion coefficients} $\gamma_{%
%TCIMACRO{\TeXButton{boldmu}{\mbox{\boldmath$ \mu$}}}%
%BeginExpansion
\mbox{\boldmath$ \mu$}%
%EndExpansion
,%
%TCIMACRO{\TeXButton{boldnu}{\mbox{\boldmath$ \nu$}}}%
%BeginExpansion
\mbox{\boldmath$ \nu$}%
%EndExpansion
,c\left(  b\right)  }^{b}\left(  \zeta\right)  $. We obtain the bound (cf.
(\ref{defnt}))%
\[
\sum_{b=\left(  t,s\right)  \in\mathcal{P}_{\operatorname*{far}}}\left(
\sharp\mathbb{N}_{t}\right)  \left(  \sharp\mathbb{N}_{s}\right)  =\sum
_{\ell=0}^{L}k_{\ell}^{2}\sum_{t\in\mathcal{T}_{\ell}}\sharp\mathcal{P}%
_{\operatorname*{right}}^{\operatorname*{far}}\left(  t\right)  \overset
{\text{Lem. \ref{lemma:inadmin}, (\ref{estRl2})}}{\leq}\hat{C}%
_{\operatorname*{sp}}\sum_{\ell=0}^{L}k_{\ell}^{2}\sum_{t\in\mathcal{T}_{\ell
}}R_{\ell}^{2}\leq k_{L}^{2}\hat{C}_{\operatorname*{sp}}\sum_{t\in
\mathcal{T}_{\mathcal{I}}}R_{t}^{2}.
\]
We may argue as in Part 1 to get the assertion.%
%TCIMACRO{\TeXButton{End Proof}{\endproof}}%
%BeginExpansion
\endproof
%EndExpansion

\begin{lemma}
[Nearfield matrix]\label{LemNear}Let the set of directions be constructed
according to \cite[Rem. 3]{borm2018hybrid}. Under assumptions (\ref{cpconst}),
(\ref{Cwidth}), (\ref{10}), (\ref{12}-\ref{15b}), (\ref{delta_ell_cond}) and
assuming that $L:=\operatorname*{depth}\left(  \mathcal{T}_{\mathcal{I}%
}\right)  >0$, there exists a constant $C_{\operatorname*{near}}$ such that
the number of non-zero nearfield matrix entries is bounded from above by%
\[
C_{\operatorname*{near}}\left(  k_{L}+1\right)  \left(  \frac{n}{\eta_{2}^{2}%
}+\min\left\{  \left(  1+L\right)  \left(  k_{L}+1\right)  \left(
\frac{\left\vert \operatorname{Im}\zeta\right\vert }{\eta_{2}}\right)
^{2},\frac{n}{\eta_{3}}\frac{\left\vert \operatorname{Im}\zeta\right\vert
}{\operatorname{Re}\zeta}\right\}  \right)  .
\]

\end{lemma}

%

%TCIMACRO{\TeXButton{Proof}{\proof}}%
%BeginExpansion
\proof
%EndExpansion
Let $\left(  t,s\right)  \in\mathcal{P}_{\operatorname*{near}}$. This implies
that $t$ or $s$ belongs to $\mathcal{L}_{\mathcal{I}}$ and we assume that
$t\in\mathcal{L}_{\mathcal{I}}\cap\mathcal{T}_{\ell}$ for some $0\leq\ell\leq
L$. From (\ref{15b}) we know that $\sharp\hat{t}\leq C_{\operatorname*{rs}%
}\left(  k_{L}+1\right)  $. The construction of the block partition implies
that also $s\in\mathcal{T}_{\ell}$. Then, we conclude as in the second
inequality of (\ref{2ndineququ})%
\begin{align}
\sharp\hat{s}  &  \leq\frac{\left\vert \omega_{s}\right\vert \left(
C_{\operatorname*{loc},p}+1\right)  }{c_{\operatorname*{sr}}h_{\min}^{2}%
}\overset{\text{(\ref{Cwidth})}}{\leq}\frac{\operatorname*{diam}%
\nolimits^{2}B_{\ell}\left(  C_{\operatorname*{loc},p}+1\right)
}{c_{\operatorname*{vol}}c_{\operatorname*{sr}}h_{\min}^{2}}\leq
\frac{C_{\operatorname*{vol}}\left\vert w_{t}\right\vert \left(
C_{\operatorname*{loc},p}+1\right)  }{c_{\operatorname*{vol}}%
c_{\operatorname*{sr}}h_{\min}^{2}}\label{estshat}\\
&  \leq\frac{C_{\operatorname*{vol}}\left(  C_{\operatorname*{loc}%
,p}+1\right)  }{c_{\operatorname*{vol}}c_{\operatorname*{sr}}}\frac
{h_{\mathcal{G}}^{2}}{h_{\min}^{2}}\sharp\hat{t}\overset{\text{(\ref{quconst}%
)}}{\leq}C_{\operatorname{level}}\sharp\hat{t}\leq C_{\operatorname{level}%
}C_{\operatorname*{rs}}\left(  k_{L}+1\right) \nonumber
\end{align}
with\ $C_{\operatorname{level}}:=\frac{C_{\operatorname*{qu}}^{2}%
C_{\operatorname*{vol}}\left(  C_{\operatorname*{loc},p}+1\right)
}{c_{\operatorname*{vol}}c_{\operatorname*{sr}}}$. This allows to estimate the
number of non-zero entries in the nearfield matrix by%
\begin{align*}
&  \sharp\left\{  \left(  i,j\right)  \in\mathcal{I}\mid K_{i,j}%
^{\operatorname*{near}}\neq0\right\}  \leq\sum_{t\in\mathcal{L}_{\mathcal{I}}%
}\sum_{s\in\mathcal{P}_{\operatorname*{right}}^{\operatorname*{near}}\left(
t\right)  \cup\mathcal{P}_{\operatorname*{left}}^{\operatorname*{near}}\left(
t\right)  }\left(  \sharp\hat{s}\right)  \left(  \sharp\hat{t}\right) \\
&  \qquad\overset{\text{Lem. \ref{lemma:basis}}}{\leq}C_{\operatorname*{near}%
}\left(  k_{L}+1\right)  \left(  \frac{n}{\eta_{2}^{2}}+\min\left\{  \left(
1+L\right)  \left(  k_{L}+1\right)  \left(  \frac{\left\vert \operatorname{Im}%
\zeta\right\vert }{\eta_{2}}\right)  ^{2},\frac{n}{\eta_{3}}\frac{\left\vert
\operatorname{Im}\zeta\right\vert }{\operatorname{Re}\zeta}\right\}  \right)
\end{align*}
for $C_{\operatorname*{near}}:=2C_{\sharp}C_{\operatorname{level}%
}C_{\operatorname*{rs}}^{2}$.
%TCIMACRO{\TeXButton{End Proof}{\endproof}}%
%BeginExpansion
\endproof
%EndExpansion

\begin{corollary}
[Nearfield matrix with resolution condition]\label{CorNMRC}Let the assumptions
of Lem\-ma \ref{LemNear} be satisfied and assume the \emph{resolution
condition}%
\[
\left\vert \operatorname{Im}\zeta\right\vert \operatorname*{diam}B_{t}\leq
C_{\operatorname*{res}}\quad\forall t\in\mathcal{L}_{\mathcal{I}}.
\]
Then, the number of non-zero nearfield matrix entries is bounded from above by%
\[
C_{\operatorname{level}}C_{\operatorname*{rs}}^{2}C_{\operatorname*{near}%
}^{\sharp}\left(  k_{L}+1\right)  n.
\]

\end{corollary}

%

%TCIMACRO{\TeXButton{Proof}{\proof}}%
%BeginExpansion
\proof
%EndExpansion
We employ Corollary \ref{LemSparsRes} to estimate the number of non-zero
nearfield matrix entries from above by%
\begin{align*}
&  C_{\operatorname*{rs}}\left(  k_{L}+1\right)  C_{\operatorname*{near}%
}^{\sharp}\frac{n}{k_{L}+1}\max_{t\in\mathcal{L}_{\mathcal{I}}}\left\{
\sharp\hat{s}\mid s\in\left(  \mathcal{P}_{\operatorname*{left}}%
^{\operatorname*{near}}\left(  t\right)  \cup\mathcal{P}%
_{\operatorname*{right}}^{\operatorname*{near}}\left(  t\right)  \right)
\right\} \\
&  \overset{\text{(\ref{estshat})}}{\leq}C_{\operatorname*{rs}}^{2}%
C_{\operatorname*{near}}^{\sharp}C_{\operatorname{level}}\left(
k_{L}+1\right)  n.
\end{align*}%
%TCIMACRO{\TeXButton{End Proof}{\endproof}}%
%BeginExpansion
\endproof
%EndExpansion

\subsection{Computational Complexity}

In Section \ref{SubSecStorage} we have described the quantities which have to
be computed and stored for the directional $\mathcal{H}^{2}$ matrix
representation for the acoustic single layer operator with complex frequency
and estimated their cardinalities. For the computational complexity the effort
for generating these quantities has to be taken into account.

\begin{remark}
[Transfer matrices]Since each cluster in $\mathcal{T}_{\mathcal{I}}%
\backslash\mathcal{L}_{\mathcal{I}}$ has a least two sons, it follows by a
geometric sum argument that $\sharp\mathcal{T}_{\mathcal{I}}\leq2\left(
\sharp\mathcal{L}_{\mathcal{I}}\right)  $. Hence, the number of the
$q_{t^{\prime},%
%TCIMACRO{\TeXButton{boldmu}{\mbox{\boldmath$ \mu$}}}%
%BeginExpansion
\mbox{\boldmath$ \mu$}%
%EndExpansion
,%
%TCIMACRO{\TeXButton{boldnu}{\mbox{\boldmath$ \nu$}}}%
%BeginExpansion
\mbox{\boldmath$ \nu$}%
%EndExpansion
}$ in (\ref{defq}) can be bounded from above by%
\[
\left(  \sharp\mathcal{T}_{\mathcal{I}}\right)  \left(  k_{L}+1\right)
^{2}\leq2\left(  \sharp\mathcal{L}_{\mathcal{I}}\right)  \left(
k_{L}+1\right)  ^{2}\overset{\text{(\ref{estleaves})}}{\leq}%
2C_{\operatorname*{rs}}n\left(  k_{L}+1\right)  .
\]
It requires the evaluation of the tensorized \v{C}eby\v{s}ev polynomials for
the clusters $t$ at the \v{C}eby\v{s}ev nodes of their sons. For the efficient
evaluation of \v{C}eby\v{s}ev polynomials we refer, e.g., to
\cite{Tref_Approx}, and denote the computational complexity for computing all
transfer matrices by%
\[
2C_{\operatorname*{Cheby}}C_{\operatorname*{rs}}n\left(  k_{L}+1\right)  ,
\]
where $C_{\operatorname*{Cheby}}$ depends algebraically on the expansion order
while an explicit estimate depends on the chosen evaluation method. We do not
elaborate on this issue here but refer to \cite{Tref_Approx} instead.
\end{remark}

\begin{remark}
[Expansion coefficients]\label{RemExpCoeff}The number of the expansion
coefficients\textbf{ }$\gamma_{%
%TCIMACRO{\TeXButton{boldmu}{\mbox{\boldmath$ \mu$}}}%
%BeginExpansion
\mbox{\boldmath$ \mu$}%
%EndExpansion
,%
%TCIMACRO{\TeXButton{boldnu}{\mbox{\boldmath$ \nu$}}}%
%BeginExpansion
\mbox{\boldmath$ \nu$}%
%EndExpansion
,c}^{b}\left(  \zeta\right)  $ can be bounded by (\ref{lemma:basisc}) and we
distinguish between three scenarios:

\begin{enumerate}
\item Sectorial case $\frac{\left\vert \operatorname{Im}\zeta\right\vert
}{\operatorname{Re}\zeta}\leq\alpha$. Then the cardinality is of order
$\mathcal{O}\left(  nk_{L}\right)  $,

\item In the high frequency case, i.e., $\operatorname{Re}\zeta=\mathcal{O}%
\left(  1\right)  $ and if the resolution condition (\ref{rescond}) is
satisfied, the cardinality is of order $\mathcal{O}\left(  \left(  1+L\right)
k_{L}^{2}n\right)  $. This shows that the asymptotic complexity of our
algorithm with the new admissibility condition is the same as the algorithm
for purely imaginary wave numbers considered in \cite[Thm. 1]{Boerm2015}.

\item In the high frequency case, i.e., $\operatorname{Re}\zeta=\mathcal{O}%
\left(  1\right)  $ and if the resolution condition (\ref{rescond}) is
violated we have $\mathcal{O}\left(  \left(  1+L\right)  k_{L}^{2}\left\vert
\operatorname*{Im}\zeta\right\vert ^{2}\right)  $.
\end{enumerate}

The evaluation of (\ref{expcoeff}) per coefficient has a computational cost of
$\mathcal{O}\left(  1\right)  $.
\end{remark}

\begin{remark}
[Basic farfield coefficients]\label{RemBFFC}The number of\textbf{ }basic
farfield coefficients $J_{%
%TCIMACRO{\TeXButton{boldmu}{\mbox{\boldmath$ \mu$}}}%
%BeginExpansion
\mbox{\boldmath$ \mu$}%
%EndExpansion
,c}^{t}$ is estimated in (\ref{lemma:basisb}). It requires the integration of
the expansion function $\Phi_{%
%TCIMACRO{\TeXButton{boldmu}{\mbox{\boldmath$ \mu$}}}%
%BeginExpansion
\mbox{\boldmath$ \mu$}%
%EndExpansion
,c}^{t}$ multiplied by basis functions. Although exact integration is feasible
on plane triangles we recommend to use tensor Gauss rules on triangles which
are transformed to squares by simplex coordinates. The order depends on the
degree $m_{t}$ of the \v{C}eby\v{s}ev polynomials on the leaves $t\in
\mathcal{L}_{\mathcal{I}}$ and the resolution condition: if the resolution
condition (\ref{rescond}) is satisfied the plane wave in the integrand is
non-oscillatory on the panels and does not cause an increase of the required
quadrature order. If the condition is violated the number of Gauss points has
to take into account the wave number. Alternatively, more specialized
quadrature methods could be employed for highly oscillatory integrals (see,
e.g., \cite{melenk_filon}). We do not discuss this issue here in detail but
assume that there exists a constant $C_{\operatorname*{qb}}$ such that the
computational complexity to compute all basis farfield coefficients is bounded
by $\mathcal{O}\left(  C_{\operatorname*{qb}}C_{\operatorname*{ffc}}%
k_{L}\left(  n+k_{L}\left(  \operatorname{Im}\zeta\right)  ^{2}\right)
\right)  $.
\end{remark}

\begin{remark}
[Nearfield matrix]The number of non-zero nearfield matrix entries is estimated
in Lemma \ref{LemNear} and Corollary \ref{CorNMRC}. Typically, numerical
quadrature is employed to approximate the integrals in (\ref{defK}) on
$\operatorname*{supp}b_{i}\times\operatorname*{supp}b_{j}$. To take into
account the singularity of the kernel functions, we recommend to use the
quadrature rules described in \cite{SauterDiss}, \cite{ErSa},
\cite{SauterSchwab2010}, where the computational effort per integral behaves
proportionally to $\left(  \log\frac{1}{\varepsilon}\right)  ^{4}$. As in
Remark \ref{RemBFFC} the number of quadrature points has to take into account
the wave number only in the case that the resolution condition (\ref{rescond})
is violated or quadrature techniques for highly oscillatory integrals should
be applied. We denote the computational complexity per non-zero nearfield
matrix entry by $C_{\operatorname*{Q}}$.

\begin{enumerate}
\item Sectorial case: The computational complexity in this case is of order
$\mathcal{O}\left(  C_{\operatorname*{Q}}nk_{L}\right)  $.

\item Non-sectorial case, resolution condition satisfied. Then, Corollary
\ref{CorNMRC} implies that the computational cost is of order $\mathcal{O}%
\left(  C_{\operatorname*{Q}}nk_{L}\right)  $.

\item Non-sectorial case, resolution condition violated. Then, the
computational cost can be estimated by using Lemma \ref{LemNear} by
$\mathcal{O}\left(  C_{\operatorname*{Q}}k_{L}^{2}\left\vert \operatorname{Im}%
\zeta\right\vert ^{2}\right)  $.
\end{enumerate}
\end{remark}

\begin{remark}
We have not taken into account the computational cost of the re-compression
algorithm because this would be a repetition of Section 5 in \cite{Boerm2015}
and Section 4 in \cite{borm2018hybrid}.
\end{remark}

\section{Numerical Experiments \label{SecNumExp}}

In order to illustrate how our theoretical results compare to practical
experiments, we consider the three-dimensional unit sphere $\Gamma
=\{x\in\mathbb{R}^{3}\ :\ \|x\|_{2}=1\}$, approximated by regularly refining
the eight triangular faces of the double pyramid $\{x\in\mathbb{R}%
^{3}\ :\ |x_{1}|+|x_{2}|+|x_{3}|=1\}$ and projecting the resulting vertices to
the sphere $\Gamma$. This yields a surface mesh with $n\in\mathbb{N}$ triangles.

We approximate the single-layer potential matrix%
\[
K_{i,j}\left(  \zeta\right)  :=\left(  V\left(  \zeta\right)  b_{j}%
,b_{i}\right)
\]
for piecewise constant basis functions $(b_{i})_{i=1}^{n}$ on the surface triangles.

Our approximation scheme uses a constant number $m$ of interpolation points
per coordinate for a total of $m^{3}$ points for a bounding box, and the
admissibility conditions (\ref{adm}) with $\eta_{1}=10$, $\eta_{2}=2$ and
$\eta_{3}=1/2$.

In all experiments we rely on algebraic recompression \cite{borm2018hybrid} to
reduce the storage requirements without significantly changing the
approximation error or the number of blocks.

\begin{figure}
[ptb]
\begin{center}
\includegraphics[
height=3.0364in,
width=5.0548in
]%
{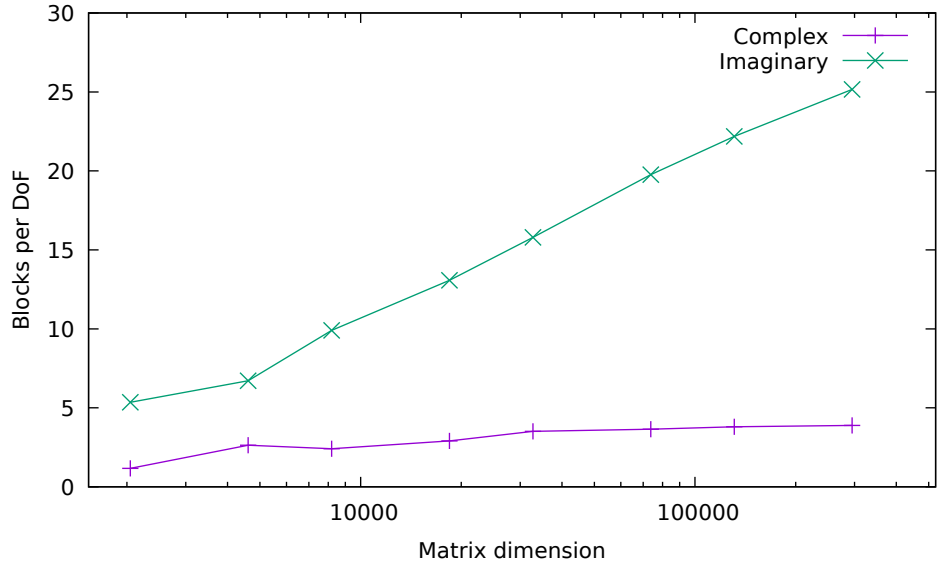}%
\caption{Number of blocks per degree of freedom depending on the matrix
dimension $n$.}%
\label{fi:blocks}%
\end{center}
\end{figure}

In a first experiment, we compare the pure Helmholtz case $\zeta
=\alpha\operatorname*{i}$ with the damped case $\zeta=\alpha+\alpha
\operatorname*{i}$, where $\alpha=\sqrt{n/128}$ guarantees $\alpha
h_{\mathcal{G}}\approx0.6$, i.e., approximately ten mesh elements per
wavelength. Figure~\ref{fi:blocks} shows the number of blocks $\#\mathcal{P}$
per degree of freedom for the purely imaginary case (labeled \textquotedblleft
Imaginary\textquotedblright) and mixed case (labeled \textquotedblleft
Complex\textquotedblright).%
%TCIMACRO{\FRAME{ftbpFU}{5.0548in}{3.0364in}{0pt}{\Qcb{Number of blocks per
%degree of freedom depending on the matrix dimension $n$.}}{\Qlb{fi:blocks}%
%}{fi_blocks.eps}{\special{ language "Scientific Word";  type "GRAPHIC";
%width 5.0548in;  height 3.0364in;  depth 0pt;  original-width 5.0004in;
%original-height 2.9922in;  cropleft "0";  croptop "1";  cropright "1";
%cropbottom "0";  filename 'fi_blocks.eps';file-properties "XNPEU";}}}%
%BeginExpansion
%EndExpansion
\begin{table}[ptb]%
\[%
\begin{array}
[c]{rr|rr|rr}%
n & \alpha & \multicolumn{2}{r|}{\zeta=\alpha+\alpha\operatorname*{i}} &
\multicolumn{2}{r}{\zeta=\alpha\operatorname*{i}}\\
&  & \#\mathcal{P} & \#\mathcal{P}/n & \#\mathcal{P} & \#\mathcal{P}/n\\\hline
2048 & 4 & 2389 & 1.17 & 10965 & 5.35\\
4608 & 6 & 12125 & 2.63 & 30941 & 6.71\\
8192 & 8 & 19733 & 2.41 & 81109 & 9.90\\
18432 & 12 & 53525 & 2.90 & 240917 & 13.07\\
32768 & 16 & 114949 & 3.51 & 517621 & 15.80\\
73728 & 24 & 268773 & 3.65 & 1457349 & 19.77\\
131072 & 32 & 498725 & 3.80 & 2908053 & 22.19\\
294912 & 48 & 1145093 & 3.88 & 7419477 & 25.16
\end{array}
\]
\end{table}

Since we are using a logarithmic scale for the matrix dimension $n$,
Figure~\ref{fi:blocks} suggests that the number of blocks grows like
$\mathcal{O}(n\log n)$ in the pure Helmholtz case, but only like
$\mathcal{O}(n)$ for the Helmholtz case with decay, in accordance with our
theoretical results.

Next we consider the dependence of the matrix approximation error, estimated
in the spectral norm by a number of steps of the power iteration for the
self-adjoint matrix $(\mathbf{K}\left(  \zeta\right)  -\widetilde
{\mathbf{K}\left(  \zeta\right)  })^{\ast}(\mathbf{K}\left(  \zeta\right)
-\widetilde{\mathbf{K}\left(  \zeta\right)  })$, on the interpolation order
$m$. Standard polynomial interpolation theory predicts that the asymptotic
\emph{rate} of convergence should be the same for all matrix dimensions $n$,
while the total error also depends on the mesh parameter.%
%TCIMACRO{\FRAME{ftbpFU}{5.0548in}{3.0364in}{0pt}{\Qcb{Convergence for
%$\zeta=\alpha+\alpha\operatorname*{i}$ and different dimensions.}%
%}{\Qlb{fi:conv_dim}}{fi_conv_dim.eps}{\special{ language "Scientific Word";
%type "GRAPHIC";  maintain-aspect-ratio TRUE;  display "USEDEF";
%valid_file "F";  width 5.0548in;  height 3.0364in;  depth 0pt;
%original-width 5.0004in;  original-height 2.9922in;  cropleft "0";
%filename 'fi_conv_dim.eps';file-properties "XNPEU";}}}%
%BeginExpansion
\begin{figure}
[ptb]
\begin{center}
\includegraphics[
height=3.0364in,
width=5.0548in
]%
{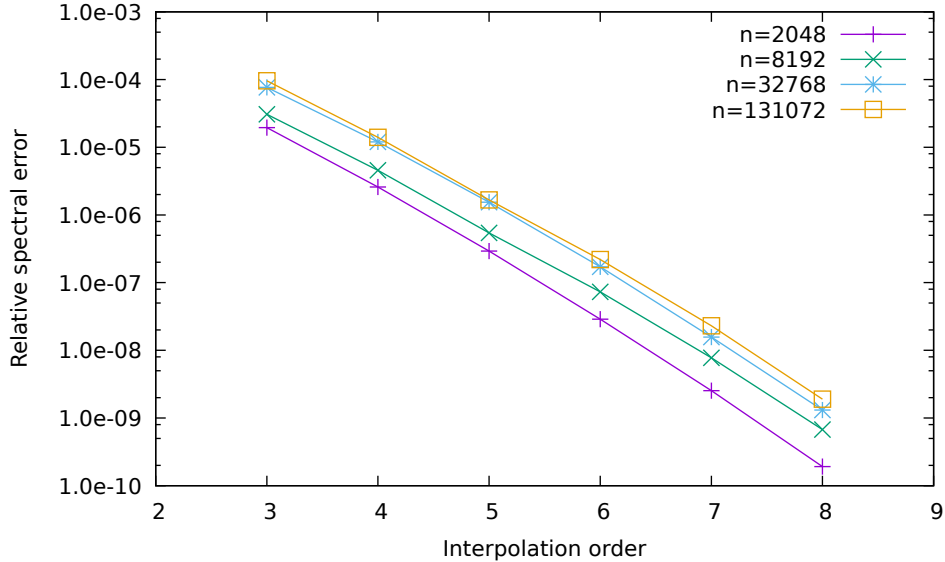}%
\caption{Convergence for $\zeta=\alpha+\alpha\operatorname*{i}$ and different
dimensions.}%
\label{fi:conv_dim}%
\end{center}
\end{figure}
%EndExpansion
Figure~\ref{fi:conv_dim} uses a logarithmic scale for the relative spectral
error on the vertical axis and a linear scale for the interpolation order on
the horizontal. We can observe the expected exponential convergence, and we
can also see that the relative error grows slowly as $n$ increases. This
latter effect can be contributed to the fact that our error estimate contains
the factor $1/\mbox{dist}(B_{t},B_{s})$ that grows like $1/h_{\mathcal{G}}$ as
the mesh is refined.

We can conclude that in the mixed case $\zeta=\alpha+\alpha\operatorname*{i}$
the number of blocks is $\mathcal{O}(n)$, while the error shows stable
convergence. Essentially the performance of the algorithm is comparable to
standard interpolation for the Laplace kernel.

\begin{figure}
[ptb]
\begin{center}
\includegraphics[
height=3.0364in,
width=5.0548in
]%
{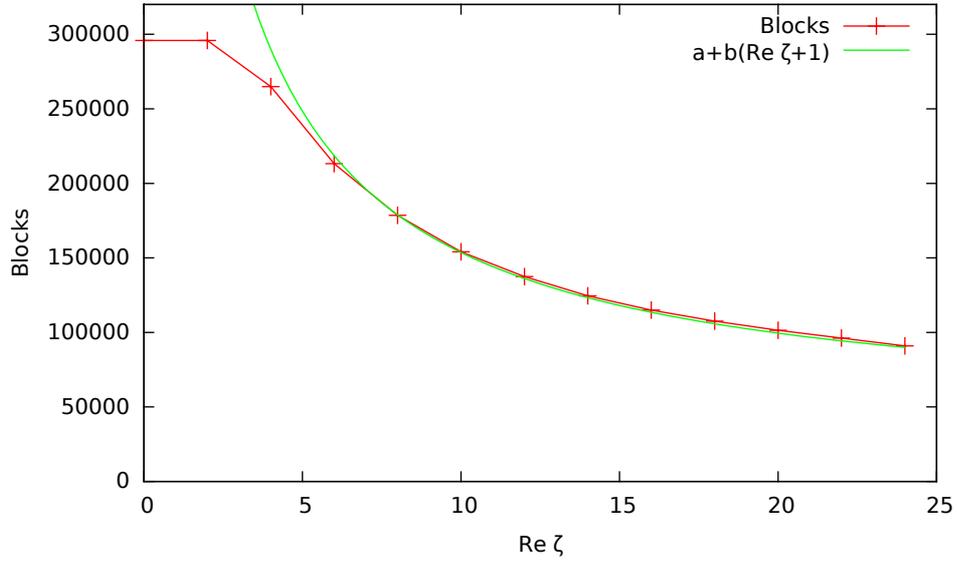}%
\caption{Number of blocks depending on the size of the real part
$\mbox{Re}\zeta$. The values of $a$ and $b$ in the fit are as in
(\ref{nuefit}): $a=40000$ and $b=1250000$.}%
\label{fi:blocks2}%
\end{center}
\end{figure}

In the context of our analysis, the dependence of the complexity and the
accuracy on the damping factor $\operatorname{Re}\zeta$ is of particular
interest. To take a closer look, we fix $n=32\,768$, $\alpha=16$, and consider
$\zeta=\nu+\alpha\operatorname*{i}$ with $\nu\in\{0,2,\ldots,24\}$. The
corresponding block numbers are shown in Figure~\ref{fi:blocks2}.%
%TCIMACRO{\FRAME{ftbpFU}{5.0548in}{3.0364in}{0pt}{\Qcb{Number of blocks
%depending on the size of the real part $\mbox{Re}\zeta$. The values of $a$ and
%$b$ in the fit are as in (\ref{nuefit}): $a=40000$ and $b=1250000$.}%
%}{\Qlb{fi:blocks2}}{fi_blocks2.eps}{\special{ language "Scientific Word";
%type "GRAPHIC";  maintain-aspect-ratio TRUE;  display "USEDEF";
%valid_file "F";  width 5.0548in;  height 3.0364in;  depth 0pt;
%original-width 5.0004in;  original-height 2.9922in;  cropleft "0";
%filename 'fi_blocks2.eps';file-properties "XNPEU";}}}%
%BeginExpansion
%EndExpansion

Based on our theoretical results we expect that the number of blocks is
proportional to $n+\min\left\{  \left\vert \operatorname*{Im}\zeta\right\vert
^{2}\log n,\frac{\left\vert \operatorname*{Im}\zeta\right\vert }%
{\operatorname{Re}\zeta}\right\}  $, and comparing the numerical results with
the function
\begin{equation}
\nu\mapsto40\,000+\frac{1\,250\,000}{\nu+1} \label{nuefit}%
\end{equation}
suggests that the prediction is quite sharp for larger values of $\nu$.

\begin{figure}
[ptb]
\begin{center}
\includegraphics[
height=3.0364in,
width=5.0548in
]%
{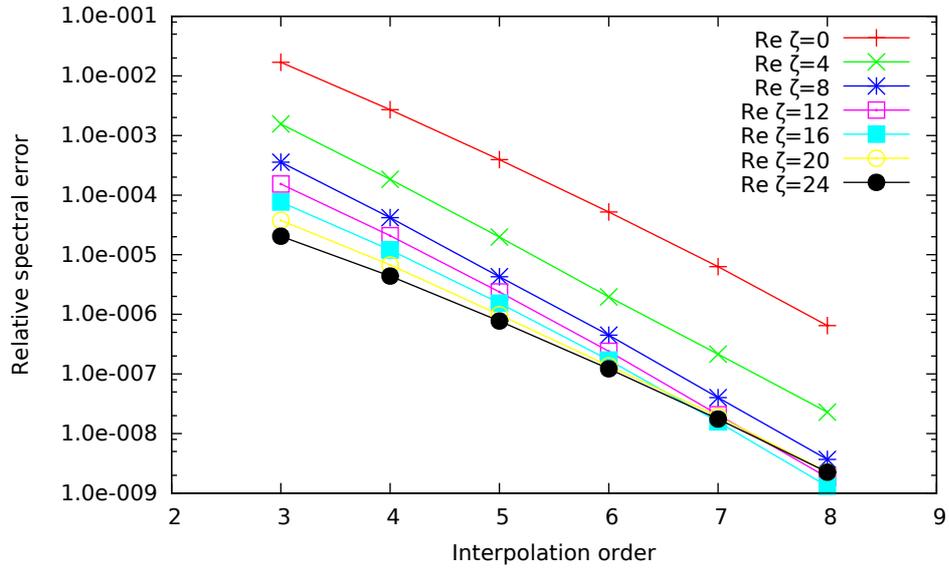}%
\caption{Spectral error versus $\operatorname{Re}\zeta$ for different
interpolation orders.}%
\label{fi:conv_nu}%
\end{center}
\end{figure}

Of course we are also interested in the dependence of the interpolation error
on the real part $\nu=\operatorname{Re}\zeta$ of $\zeta$. The relative
spectral errors for $\nu\in\{0,4,\ldots,24\}$ and interpolation orders
$m\in\{3,\ldots,8\}$ are given in Figure~\ref{fi:conv_nu}.%
%TCIMACRO{\FRAME{ftbpFU}{5.0548in}{3.0364in}{0pt}{\Qcb{Spectral error versus
%$\operatorname{Re}\zeta$ for different interpolation orders.}}%
%{\Qlb{fi:conv_nu}}{fi_conv_nu.eps}{\special{ language "Scientific Word";
%type "GRAPHIC";  maintain-aspect-ratio TRUE;  display "USEDEF";
%valid_file "F";  width 5.0548in;  height 3.0364in;  depth 0pt;
%original-width 5.0004in;  original-height 2.9922in;  cropleft "0";
%filename 'fi_conv_nu.eps';file-properties "XNPEU";}}}%
%BeginExpansion
%EndExpansion

We can see that the rates of convergence are similar for the different values
of $\nu$, while the relative errors decay as $\nu$ grows.

\bibliographystyle{abbrv}
\bibliography{acompat,mlf_stas_fast}

\end{document}